\documentclass{amsart}

\usepackage{amssymb}
\usepackage[english]{babel}
\usepackage[dvips]{graphicx}
\usepackage[lined,ruled]{algorithm2e}

\usepackage[numeric, short-months]{amsrefs}

\newtheorem{theorem}{Theorem}[section]
\newtheorem{corollary}[theorem]{Corollary}
\newtheorem{proposition}[theorem]{Proposition}
\newtheorem{lemma}[theorem]{Lemma}

\theoremstyle{definition}
\newtheorem{definition}[theorem]{Definition}
\newtheorem{problem}[theorem]{Problem}

\theoremstyle{remark}
\newtheorem{remark}[theorem]{Remark}

\numberwithin{equation}{section}

\newcommand{\Alphabet}[1]{ {X} _{#1} }
\newcommand{\InfiniteSet}[1][m]{ { \Alphabet{#1} ^ \omega } }
\newcommand{\FiniteSet}[1][m]{ { \Alphabet{#1} ^ \ast } }
\newcommand{\States}[1]{ {Q} _{#1} }

\newcommand{\AutomataSet}[2]{ \Automaton[{#1 \times #2}] }
\newcommand{\Automaton}[1][]{ { {A} _{#1} } }
\newcommand{\MAutomaton}[3]{ \left( \Alphabet{#2}, \States{#1}, \pi _{#3}, \lambda _{#3} \right) }
\newcommand{\Intermediate}[1][m]{ {J} _{#1} }
\newcommand{\Restriction}[2] { \left. {#1} \right| _{#2} }

\newcommand{\Semigroup}[1]{ {S} _{#1} }
\newcommand{\ATMonoid}[1]{ \left\langle {\begin{array}{*{20}c}
   {e, f_0, f_1} & \vline &  {#1} \\
\end{array}} \right\rangle }

\newcommand{\Sequence}[2]{ { \left\{ {#1}, m \geq {#2} \right\} } }
\newcommand{\SequenceOParam}[2]{ { \left\{ {#1[m]}, m \geq {#2} \right\} } }

\newcommand{\Growth}[1]{ \gamma _{#1} }
\newcommand{\GrowthOrder}[1]{ { \left[ {#1} \right] } }
\newcommand{\GrowthAutomaton}[1]{ \Growth{ \Automaton[#1] } }
\newcommand{\GrowthSemigroup}[1]{ \Growth{ \Semigroup{#1} } }
\newcommand{\SphericalGrowthSemigroup}[1]{ \mathord{ \buildrel{ \lower 3pt
\hbox{$\scriptscriptstyle\frown$}} \over {\Growth{}} } _{ \Semigroup{#1} } }
\newcommand{\WordGrowthSemigroup}[1]{ \delta _{ \Semigroup{#1} } }

\newcommand{\Argument}[1]{ { \left( {#1} \right) } }
\newcommand{\n}[1][n]{ \Argument{#1} }
\newcommand{\Pair}[2]{ { \Argument{#1, #2} } }
\newcommand{\Remainder}[1] { {\left[\kern-0.14em\left[ {#1}
\right]\kern-0.14em\right] } }
\newcommand{\Remainderm}[1] { {\left[\kern-0.16em\left[ {#1}
\right]\kern-0.16em\right]_m } }
\newcommand{\Dividerm}[1] { {\left[ \frac{#1}{m} \right]} }

\newcommand{\s}{ \mathsf{s} }
\newcommand{\x}{ \mathsf{x} }

\newcommand{\Natural}{ \mathbb{ N } }
\newcommand{\Integer}{ \mathbb{ Z } }
\newcommand{\PositiveInteger}{ n \in \Natural }

\begin{document}

\title[On the Mealy Automata of the Growth Order $n ^{\log n / 2 \log m}$]{On
the $3$-state Mealy Automata over an $m$-symbol Alphabet of Growth Order
$\left[ {n ^{{\log n}/{2 \log m}}} \right]$}

\author[I. Reznykov]{Illya I. Reznykov}
\address{IKC5 ltd.\\5, Krasnogvardeyskaya str., office 2\\Kyiv\\Ukraine 02094}
\curraddr{IKC5 ltd.\\5, Krasnogvardeyskaya str., office 2\\Kyiv\\Ukraine 02094}
\email{Illya.Reznykov@ikc5.com.ua}

\author[V. Sushchansky]{Vitaliy I. Sushchansky}
\address{Institute of Mathematics\\Silesian University of Technology\\ul.
Kaszubska, 23\\44-100 Gliwice\\Poland}
\curraddr{Institute of Mathematics\\Silesian University of Technology\\ul.
Kaszubska, 23\\44-100 Gliwice\\Poland}
\email{Wital.Suszczanski@polsl.pl}

\subjclass[2000]{Primary 20M35, 68Q70; Secondary 20M20, 43A35}

\date{October 10, 2005}

\keywords{Mealy automaton, growth function, intermediate growth order}

\begin{abstract}
We consider the sequence $\SequenceOParam{\Intermediate}{2}$ of the $3$-state
Mealy automata over an $m$-symbol alphabet such that the growth function of
$\Intermediate$ has the intermediate growth order $\GrowthOrder{n ^{{\log n}/{2
\log m}}}$. For each automaton $\Intermediate$ we describe the automaton
transformation monoid $\Semigroup{\Intermediate}$, defined by it, provide
generating series for the growth functions, and consider primary properties of
$\Semigroup{\Intermediate}$ and $\Intermediate$.
\end{abstract}

\maketitle

\tableofcontents

\section{Introduction} \label{sect:introduction}

\par Objects of intermediate growth attract attention of researchers,
especially after the paper of Milnor \cite{Milnor1968-Problem}, where he raised
the question on the existence of groups of intermediate growth. The first
groups of intermediate growth were constructed by Grigorchuk in 1984
\cite{Grigorchuk1984-2_groups-English} (see also
\cite{Grigorchuk1984-p_groups-English}), and the first semigroup of
intermediate growth was constructed by Belyaev, Sesekin and Trofimov in 1977
\cite{BelyaevSesekinTrofimov1977-English} (see also \cite{LavrikMannlin2001}).
As the growth of Mealy automata is close related to the growth of automaton
transformation (semi)groups, defined by them, therefore the first example of
the Mealy automaton of intermediate growth, which is called the Grigorchuk's
automaton, follows from results of \cite{Grigorchuk1984-2_groups-English}.
Various Mealy automata of intermediate growth were found in later years (see,
for ex., \citelist{\cite{FabrykowskiGupta1991}
\cite{ReznykovSushchansky2002-Intermediate-English}}). But the properties of
the growth of groups, semigroups and Mealy automata are different in kind (see,
for ex., \citelist{\cite{Shneerson2001-Intermediate}
\cite{Reznykov2004-Composite}}).

\par In \cite{Grigorchuk1989-English} Grigorchuk proves that there exists a
lacuna in intermediate growth orders of residually $p$-groups. He shows the
following result (for definitions see Section~\ref{sect:preliminaries}):
\begin{theorem}[\cite{Grigorchuk1989-English}] \label{th:lacuna_group_growth}
Let $G$ be an arbitrary finitely generated group that is residually $p$-group
for some prime $p$, and $\Growth{G}$ be the growth function of $G$. If
$\Growth{G} \prec \exp \n[\sqrt n]$, then it has polynomial growth.
\end{theorem}

\par Moreover, there exist groups of the growth order
$\GrowthOrder{\exp \n[\sqrt n]}$, that is the lower bound of intermediate
growth orders of residually $p$-groups. Indeed,

\begin{theorem}[\cite{Grigorchuk1989-English}]
For any prime $p$ there exists a finitely generated $p$-group $G$, that the
following equality holds
\[
    \Growth{G} \sim \exp \n[\sqrt n].
\]
\end{theorem}

\par On the other hand, a set of semigroup growth orders doesn't have such
lacuna. In \cite{LavrikMannlin2001} Lavrik-M{\"a}nnlin considers the growth of
two semigroups $Q$ and $S$ that were introduced in \cite{Okninski1998} and
\cite{BelyaevSesekinTrofimov1977-English}, respectively. She proves that the
growth function of the semigroup $S$ is equivalent to $\exp \n[\sqrt n]$, and
the growth function $\Growth{Q}$ of $Q$ satisfies the following equality
\[
    \Growth{Q} \sim \exp \n[\sqrt \frac{n}{\log n}],
\]
whence the growth order of $\Growth{Q}$ is strictly less than $\GrowthOrder{
\exp \n[\sqrt n]}$.

\par The Mealy automata of intermediate growth are actively studied, too. As the
group of automaton transformations defined by a Mealy automaton is residually
finite, then it follows from Theorem~\ref{th:lacuna_group_growth} that
invertible Mealy automata have a similar growth property:

\begin{theorem}[\cite{GrigorchukNekrashevichSushchansky2000-English}]
Let $\Automaton$ be an invertible Mealy automaton over the alphabet $\left\{ 0,
1, \ldots, {p - 1} \right\}$ ($p$ is a prime number), where, for any state $q$,
the output function $\lambda \Pair{\cdot}{q}$ is a power of the cyclic
permutation $\left( 0, 1, \ldots, {p - 1} \right)$. If the growth order of
$\Automaton$ is strictly less than $\GrowthOrder{\exp \n[\sqrt n]}$ then
$\Semigroup{\Automaton}$ contains a nilpotent subsemigroup of finite index, and
$\Automaton$ has polynomial growth.
\end{theorem}

\par Hence, the growth order of an arbitrary invertible automaton of
intermediate growth is greater or equal to $\GrowthOrder{\exp \n[\sqrt n]}$.
But there are no examples of invertible Mealy automata of the intermediate
growth order $\GrowthOrder{\exp \n[\sqrt n]}$.

\par Simultaneously growth of initial Mealy automata is considered,
and it produces interesting growth orders. For example, in
\cite{GrigorchukNekrashevichSushchansky2000-English} the growth function of
``the adding machine'' as the initial Mealy automaton is considered, and there
is proved that it has the logarithmic growth order $\GrowthOrder{\log _m n}$.
But the question on the existence (non-initial) Mealy automata with logarithmic
growth is still open \cite{GrigorchukNekrashevichSushchansky2000-English}.

\par There are many interesting examples of the
growth among all (invertible and non-invertible) non-initial Mealy automata.
Let us denote the set of all $n$-state Mealy automata over the $m$-symbol
alphabet by the symbol $\AutomataSet{n}{m}$. We have created the programming
system (see \cite{ReznykovSushchansky2002-Reports-English}) and have already
modelled many automata, among them all automata from the sets
$\AutomataSet{2}{2}$, $\AutomataSet{3}{2}$, $\AutomataSet{2}{3}$, and
$\AutomataSet{2}{4}$. Analyzing these data, we have found automata with new
intermediate growth orders.

\par The smallest Mealy automaton $I_2$ of intermediate growth was
found in the set $\AutomataSet{2}{2}$ (see
\cite{ReznykovSushchansky2002-Intermediate-English}). It is proved in
\cite{ReznykovSushchansky2002-Intermediate-English} that the growth order of
the growth function $\Growth{I_2}$ satisfies the following inequalities
\[
    \GrowthOrder{\exp \Argument{\sqrt[4] n}} \le \GrowthOrder{\Growth{I_2}} \le
    \GrowthOrder{\exp \Argument{\sqrt n}}.
\]
In collaboration with Bartholdi
\cite{BartholdiReznykovSushchansky2005-Intermediate} we show the sharp
asymptotic of $\Growth{I_2}$, and prove that the following equality holds
\[
    \GrowthOrder{\Growth{I_2}} = \GrowthOrder{\exp \Argument{\sqrt n}}.
\]

\par The question on the existence of Mealy automaton of intermediate
growth such that its growth function has the growth order that is less than
$\GrowthOrder{\exp \Argument{\sqrt n}}$, was raised. Basing on the results of
calculated experiments, we set up the hypothesis that intermediate growth
orders of Mealy automata fill a lacuna between polynomial and exponential
growth orders. Moreover, there exist Mealy automata with growth orders between
polynomial growth orders of integral degrees.

\par In the paper we consider the sequence $\Sequence{\Intermediate}{2}$
of the $3$-state Mealy automata over an $m$-symbol alphabet (see
Figure~\ref{fig:automaton_intermediate}) such that the growth function of
$\Intermediate$, $m \ge 2$, has the intermediate growth order $\GrowthOrder{n
^{\frac{\log n}{2 \log m}}}$. These automata substantiate the first part of our
hypothesis. Every automaton $\Intermediate$ is an example of Mealy automaton
such that the growth order of its growth function is less than
$\GrowthOrder{\exp \Argument{\sqrt n}}$. $\Intermediate[2]$ is introduced
in~\cite{Reznykov2004-Composite} in conjecture with composite growth functions.

\par The paper has the following structure. The main results are formulated in
Section~\ref{sect:main_results}, which includes three subsections. The
automaton transformation monoid $\Semigroup{\Intermediate}$, defined by
$\Intermediate$, and its relations are considered in
Subsection~\ref{subsect:main_semigroup}. The properties of the growth of
$\Semigroup{\Intermediate}$ and $\Intermediate$ are described in
Subsection~\ref{subsect:main_growth}. There are constructed the generating
series, shown sharp asymptotics, and proved interesting arithmetic properties.
Subsection~\ref{subsect:main_sequence} is devoted to the properties of
sequences, that are defined by the sequence $\Sequence{\Intermediate}{2}$.
Preliminaries are listed in Section~\ref{sect:preliminaries}. The results
listed in the subsections of Section~\ref{sect:main_results} are proved in
Section~\ref{sect:semigroup}, Section~\ref{sect:growth_functions} and
Section~\ref{sect:sequence}, respectively. Finally, in
Section~\ref{sect:final_remarks} we consider the Mealy automaton with the
``similar'' numerical properties and discuss the sequel investigations.

\section{Main results} \label{sect:main_results}

\par
\begin{figure}[t]
  \centering
  \includegraphics*{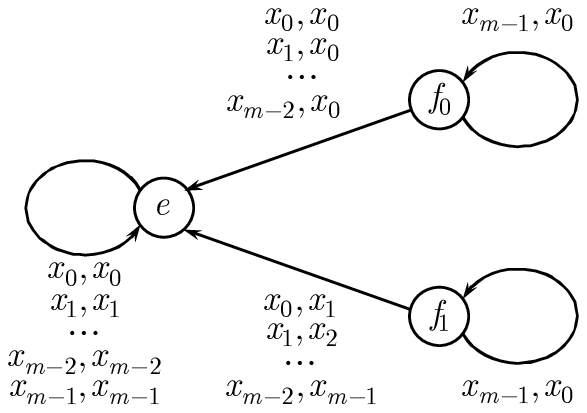}
  \caption{The automaton $\Intermediate$}
  \label{fig:automaton_intermediate}
\end{figure}

 Let $\Intermediate$, $m \ge 2$, be the $3$-state Mealy automaton over the
$m$-symbol alphabet such that its Moore diagram is shown on
Figure~\ref{fig:automaton_intermediate}. Let us denote the semigroup defined by
$\Intermediate$ by the symbol $\Semigroup{\Intermediate}$, and the growth
functions of $\Intermediate$ and $\Semigroup{\Intermediate}$ by the symbols
$\Growth{\Intermediate}$ and $\GrowthSemigroup{\Intermediate}$, respectively.

\subsection{Semigroup $\Semigroup{\Intermediate}$} \label{subsect:main_semigroup}

\par Let $m \ge 2$ be a fixed integer. The following theorem holds:

\begin{theorem} \label{th:semigroup}
The semigroup $\Semigroup{\Intermediate}$ is a monoid, and has the following
presentation:
\begin{equation*}
    \Semigroup{\Intermediate} = \ATMonoid{
        R_A \Pair{k}{p}, R_B \Argument{k}, k \ge 0, p = 1, 2, \ldots, {m -
        1}},
\end{equation*}
where the relations $R_A \Pair{k}{p}$ and $R_B \Argument{k}$ are defined by the
following equalities
\begin{multline*}
    f_0 f_1 ^{p m ^k - 1} \cdot f_0 f_1 ^ {m ^{k} - 1} f_0 \ldots  f_1 ^ {m ^2 - 1}
    f_0 f_1 ^ {m - 1} f_0\\
    = f_0 f_1 ^ {m ^{k} - 1} f_0 \ldots f_1 ^{m ^2 - 1} f_0 f_1 ^{m - 1} f_0,
\end{multline*}
and
\begin{multline*}
    f_0 f_1 ^{m ^k - 1} \cdot f_1 ^{m ^{k + 1}} f_0 f_1 ^{m ^{k} - 1} f_0
    \ldots  f_1 ^{m ^2 - 1} f_0 f_1 ^{m - 1} f_0\\
    = f_1 ^{m ^{k + 1}} f_0 f_1 ^{m ^{k} - 1} f_0 \ldots f_1 ^{m ^2 - 1} f_0
    f_1 ^{m - 1} f_0,
\end{multline*}
where $k \ge 0$, $p \ge 1$, respectively.

\par The monoid $\Semigroup{\Intermediate}$ is infinitely presented, and the
word problem may be solved in no more than quadratic time.
\end{theorem}

\begin{corollary} \label{cor:rewriting_system}
The relations
\begin{multline*}
    f_0 f_1 ^{m ^k p_{k + 2} - 1} \cdot f_1 ^{m ^{k + 1} p_{k + 1}} f_0 f_1
    ^{m ^{k} p_k - 1} f_0 f_1 ^{m ^{k - 1} p_{k - 1} - 1} f_0 \ldots f_1 ^{m
    p_1 - 1} f_0 \\
    = f_1 ^{m ^{k + 1} p_{k + 1}} f_0 f_1 ^{m ^{k} p_k - 1} f_0 f_1 ^{m ^{k -
    1} p_{k - 1} - 1} f_0 \ldots f_1 ^{m p_1 - 1} f_0,
\end{multline*}
where $k \ge 0$, $1 \le p_{k + 2} \le {m - 1}$, $p_{k + 1} \ge 0$, $p_i \ge 1$,
$i = 1, 2, \ldots, {k}$, form the rewriting system of
$\Semigroup{\Intermediate}$.
\end{corollary}

\subsection{Growth of $\Intermediate$ and $\Semigroup{\Intermediate}$}
\label{subsect:main_growth}

\par Let us denote the growth series $\sum \limits_{n \ge 0}
{\Growth{\Intermediate} \n X ^n}$ of the automaton $\Intermediate$ and the
growth series $\sum \limits_{n \ge 0} {\GrowthSemigroup{\Intermediate} \n X
^n}$ of the monoid $\Semigroup{\Intermediate}$ by the symbols $\Gamma
_{\Intermediate} \Argument{X}$ and $\Gamma _{\Semigroup{\Intermediate}}
\Argument{X}$, respectively.

\begin{theorem} \label{th:generating_functions}
The growth series $\Gamma _{\Intermediate}$ and $\Gamma
_{\Semigroup{\Intermediate}}$ coincide and admit the description
\begin{multline*}
    \Gamma _{\Intermediate} \Argument{X} = \frac{1}{\Argument{1 - X} ^2} \left(
    1 + \frac{X}{1 - X} \left( 1 + \frac{X ^m}{1 - X ^m} \left( 1 + \frac{X ^{m
    ^2}}{1 - X ^{m ^2}} \cdot \right.\right.\right.\\
    \left. \left. \left. \cdot \left( 1 + \frac{X ^{m ^3}}{1 - X ^{m ^3}}
    \left( 1 + \frac{X ^{m ^4}}{1 - X ^{m ^4}} \left(1 + \ldots \right) \right)
    \right) \right) \right) \right)
\end{multline*}

\end{theorem}

\begin{corollary} \label{cor:word_generating_function}
The word growth series $\Delta _{\Semigroup{\Intermediate}} \Argument{X} = \sum
\limits_{n \ge 0} {\WordGrowthSemigroup{\Intermediate} \n X ^n}$ of
$\Semigroup{\Intermediate}$ is defined by the following equality
\begin{multline*}
    \Delta _{\Semigroup{\Intermediate}} \Argument{X} = \frac{1}{1 - X} \left( 1
    + \frac{X}{1 - X} \left( 1 + \frac{X ^m}{1 - X ^m} \left( 1 + \frac{X ^{m
    ^2}}{1 - X ^{m ^2}} \cdot \right.\right.\right.\\
    \left. \left. \left. \cdot \left( 1 + \frac{X ^{m ^3}}{1 - X ^{m ^3}}
    \left( 1 + \frac{X ^{m ^4}}{1 - X ^{m ^4}} \left(1 + \ldots \right) \right)
    \right) \right) \right) \right)
\end{multline*}
\end{corollary}

\par Let $\Growth{}$ be an arbitrary function, and let us denote the
$i$-th finite difference of $\Growth{}$ by the symbols $\Growth{} ^{\n[i]}$, $i
\ge 1$, i.e.
\begin{align*}
    \Growth{} ^{\n[1]} \n & = \Growth{} \n - \Growth{} \n[n - 1],\\
    \Growth{} ^{\n[i]} \n & = \Growth{} ^{\n[i - 1]} \n - \Growth{}  ^{\n[i -
    1]} \n[n - 1],
\end{align*}
where $i \ge 2$, $n \ge {i + 1}$. Clearly the first difference of
$\GrowthSemigroup{\Intermediate}$ equals $\WordGrowthSemigroup{\Intermediate}$.
The arithmetic properties of $\Growth{\Intermediate}$ and
$\WordGrowthSemigroup{\Intermediate}$ are formulated in the following
corollary:

\begin{corollary} \label{cor:numerical_properties}
\begin{enumerate}
\item \label{cor_item:functional_equation} The word growth function
$\WordGrowthSemigroup{\Intermediate}$ satisfies the following \break equality
\begin{equation} \label{eq:first_difference_aymptotics}
    \WordGrowthSemigroup{\Intermediate} \n[n + 1] -
    \WordGrowthSemigroup{\Intermediate} \n =
    \WordGrowthSemigroup{\Intermediate} \n[\Dividerm{n}], \quad n \ge 0.
\end{equation}

\item \label{cor_item:gamma_one_half} The functions $\Growth{\Intermediate}$
and $\WordGrowthSemigroup{\Intermediate}$ satisfy the following equality
\begin{equation*}
    \Growth{\Intermediate} \n = \frac{1}{m} \left(
    \WordGrowthSemigroup{\Intermediate} \n[m \Argument{n + 1}] - 1 \right),
    \quad n \ge 0.
\end{equation*}

\item \label{cor_item:second_difference_partitions} Let us assume
$\Growth{\Intermediate} ^{\n[2]} \n[1] = \Growth{\Intermediate} ^{\n[2]} \n[2]
= 1$. The value $\Growth{\Intermediate} ^{\n[2]} \n$, $n \ge 1$, is equal to
the number of partitions of $n$ into ``sequential'' powers of $m$, i.e. to the
cardinality of the following set
\[
    \left\{
        \begin{array}{*{20}c}
           {p_0, p_1, \ldots, p_k} & \vline & {k \ge 0, \sum \limits_{i = 0}
           ^{k} p_i m^i = n, p_i \ge 1, i = 0, 1, \ldots, k} \\
        \end{array}
    \right\}.
\]
\end{enumerate}
\end{corollary}

\par The following theorem and corollary describe the asymptotics and the growth
orders of the functions $\Growth{\Intermediate}$ and
$\GrowthSemigroup{\Intermediate}$.

\begin{theorem} \label{th:estimates}
The growth functions have the following sharp estimates:
\begin{align*}
    \WordGrowthSemigroup{\Intermediate} \n & \sim n ^{\frac{\log {n}}{2 \log
    m}};\\
    \Growth{\Intermediate} \n  = \GrowthSemigroup{\Intermediate} \n & \sim
    \frac{1}{m} \left( m \Argument{n + 1} \right) ^{\frac{\log \left( m
    \Argument{n + 1} \right) }{2 \log m}}.
\end{align*}
\end{theorem}

\begin{corollary} \label{cor:growth_orders}
The growth orders of $\Growth{\Intermediate}$ and
$\GrowthSemigroup{\Intermediate}$ coincide, and are equal to
\[
    \GrowthOrder{\Growth{\Intermediate}} =
    \GrowthOrder{\GrowthSemigroup{\Intermediate}} = \GrowthOrder{n ^
    \frac{\log{n}}{2 \log m}}.
\]
\end{corollary}

\begin{figure}[t]
  \centering
  \includegraphics*{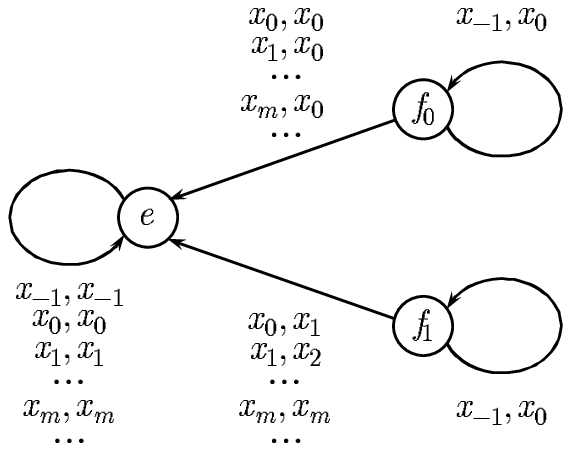}
  \caption{The automaton $\Intermediate[]'$}
  \label{fig:automaton_intermediate_infty}
\end{figure}

\subsection{The properties of $\SequenceOParam{\Intermediate}{2}$}
\label{subsect:main_sequence}

\par The sequence $\SequenceOParam{\Intermediate}{2}$ arrive in natural way
at three sequences: of the growth functions
$\Sequence{\Growth{\Intermediate}}{2}$, of the growth orders
$\Sequence{\GrowthOrder{\Growth{\Intermediate}}}{2}$, and of the automaton
transformation semigroups $\Sequence{\Semigroup{\Intermediate}}{2}$. The
following theorem characterizes boundary behavior of two of these sequences.

\begin{theorem} \label{th:sequence_growth}
\begin{enumerate}
\item \label{th_item:sequence_growth_orders} The sequence of the growth orders
$\Sequence{\GrowthOrder{\Growth{\Intermediate}}}{2}$ is a decreasing monotonic
sequence.

\item \label{th_item:sequence_growth_functions} The sequence of the growth
functions $\Sequence{\Growth{\Intermediate}}{2}$ tends pointwisely to the
function ${\Argument{n + 1}\Argument{n + 2}}/{2}$ at $m \to +\infty$.

\item \label{th_item:sequence_automata} Let $\Intermediate[]'$ be the automaton
shown on Figure~\ref{fig:automaton_intermediate_infty}. $\Intermediate[]'$ is
similar (in the sense of Definition~\ref{def:automata_similar}) to a pointwise
limit of the sequence $\SequenceOParam{\Intermediate}{2}$, and it defines the
monoid
\[
    \Semigroup{\Intermediate[]'} = \ATMonoid{ f_0 f_1 ^p f_0 =
    f_0, p \ge 0, \, f_0 f_1 ^p = f_0 f_1, p \ge 1}
\]
with the growth function $\GrowthSemigroup{\Intermediate[]'} \n = 3n$, $n \ge
1$.

Moreover, the growth function of a pointwise limit of automaton sequence
doesn't coincide with a pointwise limit of growth function sequence.
\end{enumerate}
\end{theorem}

The Item~\ref{th_item:sequence_automata} of this theorem follows from referee's
notes.

\section{Preliminaries} \label{sect:preliminaries}

\par By $\Natural$ we mean the set of non-negative integers
$\Natural = \left\{ 0, 1, 2, \ldots \right\}$.

\par We denote the remainder of a non-negative integer $p$ modulo $m$ by the
symbol $\Remainderm{p}$, and denote the integral part of a real number $r$ by
the symbol $\left[ r \right]$. Obviously for any positive integers $p, m$ the
following equality holds $p = m \Dividerm{p} + \Remainderm{p}$.

\subsection{Growth functions}

\par Let us consider the set of positive functions of a natural argument
$\Growth{} : \Natural \to \Natural$; in the sequel such functions are called
\emph{growth functions}. Let $\Growth{1} : \Natural \to \Natural$ and
$\Growth{2} : \Natural \to \Natural$ be arbitrary growth functions.

\begin{definition} \label{def:growth_order_no_greater}
The function $\Growth{1}$ has \textit{no greater growth order} (notation
$\Growth{1} \preceq \Growth{2}$) than the function $\Growth{2}$, if there exist
numbers $C_1, C_2, N_0 \in \Natural$ such that
\[
    \Growth{1} \n \le C_1 \Growth{2} \Argument{C_2 n}
\]
for any $n \ge N_0$.
\end{definition}

\begin{definition} \label{def:growth_order equivalent}
The growth functions $\Growth{1}$ and $\Growth{2}$ are equivalent or have
\textit{the same growth order} (notation $\Growth{1} \sim \Growth{2}$), if the
following inequalities hold:
\[
    \Growth{1} \preceq \Growth{2} \quad \text{and} \quad \Growth{2} \preceq
    \Growth{1}.
\]
\end{definition}

\begin{definition} \label{def:growth_order less}
The growth function $\Growth{1}$ has \textit{less growth order} (notation
$\Growth{1} \prec \Growth{2}$) than the function $\Growth{2}$, if $\Growth{1}
\preceq \Growth{2}$ but $\Growth{2} \nsim \Growth{1}$.
\end{definition}

\par The relation $\sim$ on the set of growth functions is an equivalence
relation. The equivalence class of the function $\Growth{}$ is called the
\textit{growth order} and is denoted by the symbol $\GrowthOrder{\Growth{}}$.
The relation $\preceq$ ($\prec$) induces an order relation, denoted $\le$
($<$), on equivalence classes. The growth order $\GrowthOrder{\Growth{}}$ is
called
\begin{enumerate}

\item \textit{exponential}, if $\GrowthOrder{\Growth{}} = \GrowthOrder{e ^n}$;

\item \textit{intermediate}, if $\GrowthOrder{n ^d} < \GrowthOrder{\Growth{}} <
\GrowthOrder{e ^n}$ for any $d > 0$;

\item \textit{polynomial}, if $\GrowthOrder{\Growth{}} = \GrowthOrder{n ^d}$
for some $d > 0$.
\end{enumerate}

\par The following proposition allows to compare growth orders.
\begin{proposition}[\cite{Babenko1986-English}]
\label{prop:the_same_growth_orders} Let $\Growth{1}, \Growth{2}$ be arbitrary
monotone non-decreasing growth functions. If there exist $h, a > 0$ and $b, c
\ge 0$ such that the following equality
\[
    \Growth{1} \n = h \Growth{2} \n[an + b] + c
\]
holds for all $n \ge N > 0$, then $\GrowthOrder{\Growth{1}} =
\GrowthOrder{\Growth{2}}$.
\end{proposition}

\subsection{Mealy automata}

\par Let $\Alphabet{m}$ be the $m$-symbol alphabet $\left\{ x_0, x_1, \dots,
x_{m - 1} \right\}$, $m \ge 2$. We denote the set of all finite words over
$\Alphabet{m}$, including the empty word $\varepsilon$, by the symbol
$\FiniteSet$, and denote the set of all infinite (to right) words by
$\InfiniteSet$.

\par Let $\Automaton = \MAutomaton{n}{m}{}$ be a \textit{non-initial Mealy
automaton}~\cite{Mealy1955} with the finite set of states $\States{n} = \left\{
f_0, f_1, \dots, f_{n - 1} \right\}$; input and output alphabets are the same
and are equal to $\Alphabet{m}$; $\pi : \Alphabet{m} \times \States{n} \to
\States{n}$ and $ \lambda : \Alphabet{m} \times \States{n} \to \Alphabet{m}$
are its transition and output functions, respectively. The function $\lambda$
can be extended in a natural way to a mapping $\lambda: \FiniteSet \times
\States{n} \to \FiniteSet$, and then correctly extended to a mapping $\lambda:
\InfiniteSet \times \States{n} \to \InfiniteSet$ (see, for
ex.,~\cite{Glushkov1961-English}).

\par An arbitrary Mealy automaton $\Automaton$ can be described by the Moore
diagram. The set of vertices coincides with the set of states. The edge from
the state $f$ to the state $g$ labelled by the label $x_i, x_j$ denotes that
$\pi \Pair{x_i}{f} = g$ and $\lambda \Pair{x_i}{f} = x_j$. If there are several
edges from $f$ to $g$ then we write a unique edge and join labels.

\begin{definition} \label{def:automatic_transformation}
For any state $f \in \States{n}$ the transformation $f _{\Automaton} :
\InfiniteSet \to \InfiniteSet$ defined by the equality
\[
    f _{\Automaton} \n[u] = \lambda \Pair{u}{f},
\]
where $u \in \InfiniteSet$, is called the \textit{automaton transformation}
defined by $\Automaton$ at the state $f$.
\end{definition}

\begin{definition}[\cite{Raney1958}]
\label{def:automatic_transformation_restriction} Let $f : \InfiniteSet \to
\InfiniteSet$ be an arbitrary automaton transformation, and $u \in \FiniteSet$.
The automaton transformation $\Restriction{f}{u} : \InfiniteSet \to
\InfiniteSet$, defined by
\[
    f \n[uw] = v \cdot \Restriction{f}{u} \n[w],
\]
where $w \in \InfiniteSet$ and $v$ is the beginning of $f \n[uw]$ of length
$\left| u \right|$, is called the \textit{restriction} of $f$ at the word $u$.
\end{definition}

\par The restrictions of the automaton transformation are characterized by the
following proposition.
\begin{proposition}[\cite{Raney1958}]
Let $f$ be an automaton transformation, defined by the automaton $\Automaton$
at the state $f$, $u \in \FiniteSet$ be an arbitrary finite word. Then the
restriction $\Restriction{f}{u}$ is equal to the transformation defined by
$\Automaton$ at the state $\pi \Pair{u}{f}$.
\end{proposition}

\par Let $f$ be an arbitrary state. Interpreting an
automaton transformation as an endomorphism of the rooted $m$-regular tree
(see, for ex., \cite{GrigorchukNekrashevichSushchansky2000-English}), the image
of the word $u = u_0 u_1 u_2 \cdots \in \InfiniteSet$ under the action of $f
_{\Automaton}$ can be written in the following way:
\[
    f _{\Automaton} \Argument{u_0 u_1 u_2 \dots} = \lambda \Pair{u_0}{f} \cdot
    g _{\Automaton} \Argument{u_1 u_2 \dots } = \sigma _{f} \n[u_0]
    \cdot g _{\Automaton} \Argument {u_1 u_2 \dots},
\]
where $g = \pi \Argument{u_0, f}$ and
\[
    \sigma _{f} = \left( {\begin{array}{*{20}c}
        {x_0} & {x_1} &  \ldots  & {x_{m - 1}}  \\
        {\lambda \Pair{x_0}{f}} & {\lambda \Pair{x_1}{f}} &  \ldots  & {\lambda
        \Pair{x_{m - 1}}{f}} \\
    \end{array} } \right)
\]
It means that $f _{\Automaton}$ acts on the first symbol of $u$ by the
transformation $\sigma _{f}$ over $\Alphabet{m}$, and acts on the remainder of
$u$ without its first symbol by the automaton transformation ${\pi
\Pair{u_0}{f}} _{\Automaton}$. Therefore the transformations defined by
$\Automaton$ have the following \textit{decomposition}:
\[
    f _{i} = \left( {\pi \Pair{x_0}{f_i}}, {\pi \Pair{x_1}{f_i}}, \ldots, {\pi
    \Pair{x_{m - 1}}{f_i}} \right) \sigma _{f_i},
\]
where $i = 0, 1, \dots, {n - 1}$. The Mealy automaton $\Automaton =
\MAutomaton{n}{m}{}$ defines the set
\[
    F_{\Automaton} = \left\{ f _{0}, f _{1}, \ldots, f _{n - 1} \right\}
\]
of automaton transformations over $\InfiniteSet$. The Mealy automaton
$\Automaton$ is called \textit{invertible} if all transformations from the set
$F _{\Automaton}$ are bijections. It is easy to show that $\Automaton$ is
invertible iff the transformation $\sigma _{f}$ is a permutation of
$\Alphabet{m}$ for each state $f \in \States{n}$.

\begin{definition}[\cite{Glushkov1961-English}]
The Mealy automata $\Automaton[i] = \MAutomaton{n_i}{m}{i}$ for $i = 1, 2$ are
called \textit{equivalent} if $F _{\Automaton[1]} = F _{\Automaton[2]} $.
\end{definition}

\begin{proposition}[\cite{Glushkov1961-English}]
Each class of equivalent Mealy automata over the alphabet $\Alphabet{m}$
contains, up to isomorphism, a unique automaton that is minimal with respect to
the number of states (such an automaton is called \textit{reduced}).
\end{proposition}
The minimal automaton can be found using the standard algorithm of
minimization.

\begin{definition} \label{def:automata_similar}
  The Mealy automata $\Automaton[i] = \MAutomaton{n}{m}{i}$ for $i = 1, 2$ are
  called \textit{similar} if there exist permutations $\xi \in Sym
  \Argument{\Alphabet{m}}$ and $\theta \in Sym \Argument{\States{n}}$ such that
\begin{align*}
    \theta \pi _1 (\x, f) & = \pi _2 ({\xi \x, \theta f}), & \xi \lambda _1
    (\x, f) & = \lambda _2 ({\xi \x,\theta f})
\end{align*}
for all $\x \in \Alphabet{m}$ and $f \in \States{n}$.
\end{definition}

\begin{definition}[\cite{Gecseg1986}] \label{def:automata_product}
For $i=1,2$ let $\Automaton[i] = \MAutomaton{n_i}{m}{i}$ be arbitrary Mealy
automata. The automaton $\Automaton = \left( \Alphabet{m}, \States{n _1} \times
\States{n _2}, \pi, \lambda \right)$ such that its transition and output
functions are defined by the following equalities
\begin{align*}
    \pi \Argument{\x, \Pair{f}{g}} & = \Pair{\pi _1 \Argument{\lambda _2
    \Argument{\x, g}, f}}{\pi _2 \Argument{\x, g}},\\
    \lambda \Argument{\x, \Pair{f}{g}} & = \lambda _1 \Argument{\lambda
    _2 \Argument{\x, g}, f},
\end{align*}
where $\x \in \Alphabet{m}$ and $\Pair{f}{g} \in \States{n _1} \times \States{n
_2}$, is called the \emph{product} of $\Automaton[1]$ and $\Automaton[2]$.
\end{definition}

\par We apply the automaton transformations in right to left order, that is for
arbitrary automaton transformations $f, g$ and for all $u \in \InfiniteSet$ the
equality $f \cdot g \n[u] = f \n[{g \n[u]}]$ holds.

\begin{proposition}[\cite{Gecseg1986}] \label{prop:automaton_multiple}
For any states $f \in \States{n _1}$ and $g \in \States{n _2}$ and an arbitrary
word $u \in \FiniteSet$ the following equality holds:
\[
    {\Pair{f}{g}} _{\Automaton[1] \times \Automaton[2]} \n[u] =
    {f}_{\Automaton[1]} \Argument{{g}_{\Automaton[2]} \n[u]}.
\]
\end{proposition}

\par It follows from Proposition~\ref{prop:automaton_multiple} that for the
transformations $f _{\Automaton[1]}$ and $g _{\Automaton[2]}$ the decomposition
of the product $\Pair{f}{g} _{\Automaton[1] \times \Automaton[2]}$ is defined
by:
\[
    \Pair{f}{g} _{\Automaton[1] \times \Automaton[2]} = f _{\Automaton[1]}
    \cdot g _{\Automaton[2]} = \left( h_0, h_1, \ldots, h_{m - 1} \right)
    \sigma _{f, \Automaton[1]} \sigma _{g, \Automaton[2]},
\]
where the transformation $h _i = {\pi_1 \Pair{\sigma _{g, \Automaton[2]}
\n[x_i]}{f}} _{\Automaton[1]} \cdot {\pi_2 \Pair{x_i}{g}} _{\Automaton[2]}$ for
$i = 0, 1, \dots, {m - 1}$.

\par The power $\Automaton ^n$ is defined for any automaton $\Automaton$ and
any positive integer $n$. Let us denote $\Automaton ^{\n}$ the minimal Mealy
automaton equivalent to $\Automaton ^n$. It follows from
Definition~\ref{def:automata_product} that $\left| {\States{\Automaton ^{\n}}}
\right| \le \left| {\States{\Automaton}} \right| ^n$. In addition, let
$\Automaton ^0$ be the $1$-state automaton over an $m$-symbol alphabet such
that $\sigma _{f_0}$ is the identical permutation if the semigroup
$\Semigroup{\Automaton}$ is a monoid; and $\Automaton ^0$ be the $0$-state
Mealy automaton otherwise.

\begin{definition}[\cite{Grigorchuk1988-English}] \label{def:growth_automaton}
The function $\GrowthAutomaton{}$ of a natural argument, defined by
\[
    \GrowthAutomaton{} \n = \left| {\States{\Automaton^{\n}}} \right|,
\]
where $\PositiveInteger$, is called the \textit{growth function} of the Mealy
automaton $\Automaton$.
\end{definition}

\par It is often convenient to encode the growth function in a generating
series:
\begin{definition}
Let $\Automaton$ be an arbitrary Mealy automaton. The \emph{growth series} of
$\Automaton$ is the formal power series
\[
    \Gamma _{\Automaton} \n[X] = \sum \limits _{n \ge 0} \GrowthAutomaton{} \n
    X ^n.
\]
\end{definition}

\subsection{Semigroups}

\par The necessary definitions concerning semigroups may be found in
\cite{Lallement1979}. Let $\Semigroup{}$ be a semigroup with the finite set of
generators $G = \left\{ s_0, s_1, \dots, s_{k - 1} \right\}$. \textit{The
length} of a semigroup element $\s$ is defined as a distance at the semigroup
graph from the identity in a natural metrics, that is
\[
    \ell \n[\s] = \min \limits_{l} \left\{
        \begin{array}{*{20}c}
            {\s = {s_{i_1}} {s_{i_2}} {s_{i_3}} \ldots {s_{i_l}}} & \vline &
            {s_{i_j} \in G, 1 \le j \le l} \\
        \end{array}
    \right\}.
\]
Obviously for any $\s \in \Semigroup{}$ the inequality $\ell \n[\s] > 0$ holds;
and let $\ell \n[e] = 0$ when $\Semigroup{}$ is a monoid. The \emph{normal
form} of a semigroup word is the equivalent semigroup word of minimal length.

\par Rewriting system for a semigroup is a set of equations (rules) of the
form $v = w$. A semigroup word is \emph{reduced} if it doesn't contain
occurrence of the left-hand side of a rule. The rewriting system is complete if
the set of reduced words is in bijection with the semigroup.

\par We will use several different growth functions of a semigroup. These
functions are close related with each other but they demonstrate different
properties in the case of semigroups.

\begin{definition} \label{def:growth_semigroup}
The function $\GrowthSemigroup{}$ of a natural argument $\PositiveInteger$
defined by
\[
    \GrowthSemigroup{} \n = \left| \left\{
        {\begin{array}{*{20}c}
            {s \in \Semigroup{}} & \vline & {\ell \n[s] \le n}  \\
        \end{array}} \right\} \right|
\]
is called the \textit{growth function of $\Semigroup{}$ relative to the system
$G$ of generators}.
\end{definition}

\begin{definition} \label{def:spherical_growth_semigroup}
The function $\SphericalGrowthSemigroup{}$ of a natural argument
$\PositiveInteger$ defined by
\[
    \SphericalGrowthSemigroup{} \n = \left| \left\{
        {\begin{array}{*{20}c}
            {s \in \Semigroup{}} & \vline & {s = s_{i_1} s_{i_2} \dots s_{i_n},
            \, s_{i_j} \in G, \, 1 \le j \le n}  \\
        \end{array}} \right\} \right|
\]
is called the \textit{spherical growth function of $\Semigroup{}$ relative to
the system $G$ of generators}.
\end{definition}

\begin{definition} \label{def:word_growth_semigroup}
The function $\WordGrowthSemigroup{}$ of a natural argument
$\PositiveInteger$ defined by
\[
    \WordGrowthSemigroup{} \n = \left| \left\{
        {\begin{array}{*{20}c}
            {s \in \Semigroup{}} & \vline & {\ell \n[s] = n}  \\
        \end{array}} \right\} \right|
\]
is called the \textit{word growth function of $\Semigroup{}$ relative to the
system $G$ of generators}.
\end{definition}

\par The following proposition is well-known (see,
for ex.,~\cites{Nathanson1999}):

\begin{proposition} \label{prop:equivalent_growth_semigroup}
Let $\Semigroup{}$ be an arbitrary finitely generated semigroup, and let $G_1$
and $G_2$ be systems of generators of $\Semigroup{}$.  Let us denote the growth
function of $\Semigroup{}$ relative to the set $G_i$ of generators by the
symbol $\GrowthSemigroup{i}$, for $i = 1, 2$. Then
$\GrowthOrder{\GrowthSemigroup{1}} = \GrowthOrder{\GrowthSemigroup{2}}$.
\end{proposition}

\par From
Definitions~\ref{def:growth_semigroup},~\ref{def:spherical_growth_semigroup}
and~\ref{def:word_growth_semigroup} follows that the inequalities hold
\begin{equation} \label{eq:growths_semigroup}
    \WordGrowthSemigroup{} \n \le \SphericalGrowthSemigroup{} \n \le
    \GrowthSemigroup{} \n = \sum \limits _{i = 0} ^n {\WordGrowthSemigroup{}
    \Argument{i}}, \quad \PositiveInteger.
\end{equation}

\begin{proposition} \label{prop:monoid_growth_order}
Let $\Semigroup{}$ be an arbitrary finitely generated monoid. Then
\[
    \GrowthOrder{\WordGrowthSemigroup{}} \le
    \GrowthOrder{\SphericalGrowthSemigroup{}} =
    \GrowthOrder{\GrowthSemigroup{}}.
\]
If the system $G$ of generators includes the identity, then for all $n \in
\Natural$ the equality
\[
    \SphericalGrowthSemigroup{} \n = \GrowthSemigroup{} \n
\]
holds, where the growth functions are considered relatively to the set $G$.
\end{proposition}

\par The growth function of a semigroup can be encode in a
generating series, too:
\begin{definition}
Let $\Semigroup{}$ be a semigroup generated by a finite set $G$. The
\emph{growth series} of $\Semigroup{}$ is the formal power series
\[
    \Gamma _{\Semigroup{}} \n[X] = \sum \limits _{n \ge 0}
    \GrowthSemigroup{} \n X ^n.
\]
\end{definition}
The power series $\Delta _{\Semigroup{}} \n[X] = \sum \limits_{n \ge 0}
\WordGrowthSemigroup{} \n X ^n$ can also be introduced; we then have $\Delta
_{\Semigroup{}} \n[X] = \Argument{1 - X} \Gamma _{\Semigroup{}} \n[X]$. The
series $\Delta _{\Semigroup{}}$ is called the \textit{word growth series} of
the semigroup $\Semigroup{}$.

\begin{definition} \label{def:transformation_semigroup}
Let $\Automaton = \MAutomaton{n}{m}{}$ be a Mealy automaton. A semigroup
\[
    \Semigroup{\Automaton} = \mathop{sg} \Argument{ f_{0}, f_{1}, \ldots , f
    _{n - 1}}
\]
is called the \textit{automaton transformation semigroup defined by
$\Automaton$}.
\end{definition}

\par Let $\Automaton$ be a Mealy automaton, let $\Semigroup{\Automaton}$ be the
semigroup defined by $\Automaton$, and let us denote the growth function and
the spherical growth function of $\Semigroup{\Automaton}$ by the symbols
$\GrowthSemigroup{\Automaton}$ and $\SphericalGrowthSemigroup{\Automaton}$,
respectively. From Definition~\ref{def:transformation_semigroup} we have

\begin{proposition}[\cite{Grigorchuk1988-English}]
\label{prop:equivalent_growth_functions} For any $\PositiveInteger$ the value
$\GrowthAutomaton{} \n$ is equal to the number of those elements of
$\Semigroup{\Automaton}$ that can be presented as a product of length $n$ in
the generators $\left\{ f_{0}, f_{1}, \dots, f_{n - 1} \right\}$, i.e.
\[
    \GrowthAutomaton{} \n = \SphericalGrowthSemigroup{\Automaton} \n, \,
    \PositiveInteger.
\]
\end{proposition}

\begin{proposition}
    Let $\Automaton[i]$, $i = 1, 2$, be arbitrary similar automata. Then
    $\Semigroup{\Automaton[1]}$ and $\Semigroup{\Automaton[2]}$ are isomorphic
    semigroups, and $\GrowthAutomaton{1} \n = \GrowthAutomaton{2} \n$ for all $n
    \ge 0$.
\end{proposition}

\par From this proposition and~\eqref{eq:growths_semigroup} follows that
$\GrowthAutomaton{} \n \le \GrowthSemigroup{\Automaton} \n$ for any
$\PositiveInteger$. Moreover, Mealy automata of polynomial growth such that the
equality $\GrowthOrder{\GrowthAutomaton{}} <
\GrowthOrder{\GrowthSemigroup{\Automaton}}$ holds are considered in
\cite{Reznykov2003-Polynomial}.

\section{Semigroup $\Semigroup{\Intermediate}$}
\label{sect:semigroup}

\par Let us fix $m \ge 2$ in this section.

\subsection{Semigroup relations} \label{subsect:semigroup_relations}

\par Let $\alpha_i : \Alphabet{m} \to \Alphabet{m}$, $i = 0, 1, \ldots, {m -
1}$, be the transformation such that $\alpha_i \n[\x] = x_i$ for all $\x \in
\Alphabet{m}$. Let $\sigma : \Alphabet{m} \to \Alphabet{m}$ be the permutation
such that $\sigma \n[x_i] = x_{\n[i + 1] \mod m}$ for all $i = 0, 1, \ldots, {m
- 1}$. Then $\alpha_i$ and $\sigma$ are defined by the following equalities
\begin{align*}
    \alpha _i & = \left( {\begin{array}{*{20}c}
        {x_0} & {x_1} &  \ldots  & {x_{m - 1}}  \\
        {x_i} & {x_i} &  \ldots  & {x_i}  \\
    \end{array} } \right), &
    \sigma & = \left( {\begin{array}{*{20}c}
        {x_0} & {x_1} &  \ldots  & {x_{m - 2}} & {x_{m - 1}} \\
        {x_1} & {x_2} &  \ldots  & {x_{m - 1}} & {x_0} \\
    \end{array} } \right).
\end{align*}
Using these equalities, the power of $\sigma$ is defined by the following
equality
\[
    \sigma ^i = \left( {\begin{array}{*{20}c}
        {x_0} & {x_1} &  \ldots  & {x_{m - 2}} & {x_{m - 1}} \\
        {x_{\Remainderm{i}}} & {x_{\Remainderm{i + 1}}} &  \ldots  &
        {x_{\Remainderm{i + m - 2}}} & {x_{\Remainderm{i + m - 1}}} \\
    \end{array} } \right)
\]
for all $i \ge 0$. In addition, $\sigma ^i = \sigma ^j$ if and only if $i
\equiv j \mod m$.

\par The automaton $\Intermediate$ obviously defines the identical automaton
transformation at the state $e$, and therefore $\Semigroup{\Intermediate}$ is a
monoid. In the sequel, we assume $f ^0 = e$ for an arbitrary automaton
transformation $f$. Using these agreements, the decompositions of the
transformations $f_0$ and $f_1$ are defined by the following equalities
\begin{equation} \label{eq:unrolled_form}
\begin{aligned}
    f_0 &= \left( e, e, \ldots, e, f_0 \right) \alpha _0,\\
    f_1 &= \left( e, e, \ldots, e, f_1 \right) \sigma.
\end{aligned}
\end{equation}

\par Let $\Integer _m = \left\{ 0, 1, \ldots, m - 1 \right\}$ and let $\eta :
\Integer _m \to \Alphabet{m}$ be a natural bijection such that $\eta \n[i] =
x_i$. The function $\eta$ can be extended to a mapping of $\Integer$ into the
set of infinite words, where each integer is considered as an $m$-adic number
written from left-to-right order and supplemented with infinite sequence of $0$
or $1$ depending on a sign.

\par It follows from~\eqref{eq:unrolled_form} that the action of $f_1$ can be
interpreted as the adding one to the input number. Namely, for any $p_0 \ge 0$
and $p_1$ we have
\begin{equation*}
    f_1 ^{p_0} \n[ {\eta \n[p_1]} ] = \eta \n[{p_0 + p_1}].
\end{equation*}
The action of the automaton transformation $f_0$ can be described in the
following way. It follows from the Moore diagram of $\Intermediate$ that $f_0$
replaces each symbol $x_{m - 1}$ till the first symbol $y \neq x_{m - 1}$ by
$x_0$, and then replaces $y$ by $x_0$. Let $p$ and $q$ be arbitrary $m$-adic
numbers:
\begin{align*}
    p & = \sum \limits _{n \ge 0} {p_n m ^n}, & q & = \sum \limits _{n \ge 0}
    {q_n m ^n},
\end{align*}
where $p_n, q_n \in \left\{ 0, 1, \ldots, {m - 1} \right\}$. Let $\& _m$ be a
binary operation such that
\begin{equation*}
    p \mathop{\& _m} q = \sum \limits _{n \ge 0} {\Argument{p_n \cdot \delta
    _{p_n q_n}} m ^n},
\end{equation*}
where $\delta _{p_n q_n}$ is a Kronecker symbol, $\delta _{p_n q_n} = 1$ if
$p_n = q_n$, and $\delta _{p_n q_n} = 0$ otherwise. Note that the operation $\&
_2$ coincides with the bitwise ``and'' operation. Then for any $p$ the
following equality holds
\begin{equation*}
    f_0 \n[ {\eta \n[p]} ] = \eta \n[{p \mathop{\& _m} \left( p + 1 \right)}].
\end{equation*}

\par The simple properties of $f_0$ and $f_1$ are described in the following
lemmas.
\begin{lemma}
    The relation $f_0 ^2 = f_0$ holds in $\Semigroup{\Intermediate}$.
\end{lemma}

\begin{lemma}
    The transformation $f_1$ is a bijection.
\end{lemma}

\begin{lemma} \label{lem:power_f_1}
For any $p \ge 0$ the following equality holds
\begin{equation*}
    f_1 ^{p} = \left( f_1 ^{\Dividerm{p}}, f_1 ^{\Dividerm{p + 1}}, \ldots, f_1
    ^{\Dividerm{p + m - 2}}, f_1 ^{\Dividerm{p + m - 1}} \right) \sigma ^p.
\end{equation*}
\end{lemma}

\begin{proof}
\par Let us prove Lemma~\ref{lem:power_f_1} by induction on $p$. For $p = 0$ we
have
\[
    f_1 ^{0} = e = \left( f_1 ^{\Dividerm{0}}, f_1 ^{\Dividerm{1}}, \ldots, f_1
    ^{\Dividerm{m - 1}} \right) \sigma ^0,
\]
and for $p > 1$ the equality follows from~\eqref{eq:unrolled_form}
\begin{align*}
    f_1 ^{p} & = \left( f_1 ^{\Dividerm{p - 1}}, f_1 ^{\Dividerm{p}}, \ldots,
    f_1 ^{\Dividerm{p + m - 2}} \right) \sigma ^{p - 1} \cdot \left( e, e,
    \ldots, e, f_1 \right) \sigma\\
    & = \left( f_1 ^{\Dividerm{p}}, f_1 ^{\Dividerm{p + 1}}, \ldots, f_1
    ^{\Dividerm{p + m - 2}}, f_1
    ^{\Dividerm{p - 1} + 1} \right) \sigma ^p. \qedhere
\end{align*}
\end{proof}

\par Using Lemma~\ref{lem:power_f_1} and~\eqref{eq:unrolled_form}, for any $i
\ge 1$ and $p \ge 1$ we have
\begin{subequations}
\begin{equation*}
    f_1 ^{p m ^i - 1} = \left( f_1 ^{p m ^{i - 1} - 1}, f_1 ^{p m ^{i - 1}},
    f_1 ^{p m ^{i - 1}}, \ldots, f_1 ^{p m ^{i - 1}} \right) \sigma ^{m - 1},
\end{equation*}
whence
\begin{align}
\notag
    f_1 ^{p m ^i - 1} f_0 & = \left( f_1 ^{p m ^{i - 1} - 1}, f_1 ^{p m ^{i - 1}
    - 1}, \ldots, f_1 ^{p m ^{i - 1} - 1}, f_1 ^{p m ^{i - 1} - 1} f_0 \right)
    \alpha _{m - 1},\\
\intertext{and}
\label{eq:unrolled_form_p_01}
    f_0 f_1 ^{p m ^i - 1} & = \left( f_0 f_1 ^{p m ^{i - 1} - 1}, f_1 ^{p m ^{i -
    1}}, f_1 ^{p m ^{i - 1}}, \ldots, f_1 ^{p m ^{i - 1}} \right) \alpha _0.
\end{align}
Let us denote
\begin{align*}
    v_k & = f_0 f_1 ^{m ^k - 1} f_0 f_1 ^{m ^{k - 1} - 1} \ldots f_0 f_1 ^{m -
    1} f_0,
\end{align*}
where $k \ge 0$ and $v_0 = f_0$. It follows from~\eqref{eq:unrolled_form_p_01}
that the transformation $v_k$ for $k \ge 1$ has the following decomposition
\begin{equation*}
\begin{aligned}
    f_0 f_1 ^{m ^k - 1} f_0 & f_1 ^{m ^{k - 1} - 1} \ldots f_0 f_1 ^{m ^2 - 1}
    f_0 f_1 ^{m - 1} f_0 =\\
    = & \left( f_0 f_1 ^{m ^{k - 1} - 1} f_0 f_1 ^{m ^{k - 2} - 1} \ldots f_0
    f_1 ^{m - 1} \cdot f_0 e \cdot e, \right.\\
    & \quad \ldots\\
    & \quad f_0 f_1 ^{m ^{k - 1} - 1} f_0 f_1 ^{m ^{k - 2} - 1} \ldots f_0
    f_1 ^{m - 1} \cdot f_0 e \cdot e,\\
    & \quad \left. f_0 f_1 ^{m ^{k - 1} - 1} f_0 f_1 ^{m ^{k - 2} - 1} \ldots f_0
    f_1 ^{m - 1} \cdot f_0 e \cdot f_0 \right) \alpha _0,
\end{aligned}
\end{equation*}
whence
\begin{equation} \label{eq:unrolled_form_v_k}
    v_k = \left( v_{k - 1}, v_{k - 1}, \ldots, v_{k - 1} \right) \alpha _0, \quad
    k \ge 1.
\end{equation}
\end{subequations}

\par Now we construct the irreducible system of semigroup relations.

\begin{proposition} \label{prop:defining_relations}
In the semigroup $\Semigroup{\Intermediate}$ the following relations hold:
\begin{multline} \label{eq:defining_relation_A}
    R_A \Pair{k}{p} : \, f_0 f_1 ^{p m ^k - 1} \cdot f_0 f_1 ^ {m ^{k} - 1} f_0
    \ldots  f_1 ^ {m ^2 - 1} f_0 f_1 ^ {m - 1} f_0 \\
    = f_0 f_1 ^ {m ^{k} - 1} f_0 \ldots f_1 ^{m ^2 - 1} f_0 f_1 ^{m - 1} f_0,
\end{multline}
and
\begin{multline} \label{eq:defining_relation_B}
    R_B \Argument{k} : \, f_0 f_1 ^{m ^k - 1} \cdot f_1 ^{m ^{k + 1}} f_0 f_1
    ^{m ^{k} - 1} f_0 \ldots  f_1 ^{m ^2 - 1} f_0 f_1 ^{m - 1} f_0 \\
    = f_1 ^{m ^{k + 1}} f_0 f_1 ^{m ^{k} - 1} f_0 \ldots f_1 ^{m ^2 - 1} f_0
    f_1 ^{m - 1} f_0,
\end{multline}
where $k \ge 0$, $p = 1, 2, \ldots, {m - 1}$.
\end{proposition}

\begin{remark}
Let us call the relations $R_A \Pair{k}{p}$ and $R_B \Argument{k}$ as the
relation of type A of length $k$ and the relation of type B of length $k$,
respectively. In addition, relations~\eqref{eq:defining_relation_A}
and~\eqref{eq:defining_relation_B} can be written in the following way
\begin{align*}
    R_A \Pair{k}{p} & :\, f_0 f_1 ^{p m ^k - 1} \cdot v_k = v_k; & R_B
    \Argument{k}&: \, f_0 f_1 ^{m ^k + m ^{k + 1} - 1} \cdot v_k = f_1 ^{m ^{k
    + 1}} \cdot v_k.
\end{align*}
\end{remark}

\begin{proof}
Let us prove the lemma by induction on $k$. For $k = 0$ the
relations~\eqref{eq:defining_relation_A} and~\eqref{eq:defining_relation_B} are
written in the following way
\begin{align*}
    R_A \Pair{0}{p} : \, f_0 f_1 ^{p - 1} f_0 & = f_0; & R_B \Argument{0} : \,
    f_0 f_1 ^{m} f_0 & = f_1 ^{m} f_0.
\end{align*}
Let $1 \le p \le {m - 1}$, and it follows from Lemma~\ref{lem:power_f_1} that
the equalities hold
\begin{multline*}
    f_0 f_1 ^{p - 1} f_0 = \left( e, \ldots, e, f_0 \right) \alpha _0 \cdot
    \left( f_1 ^{\Dividerm{p - 1}}, \ldots, f_1 ^{\Dividerm{p - 1}}, f_1
    ^{\Dividerm{p - 1}} f_0 \right) \alpha _{\Remainderm{p - 1}}\\
    = \left( e, e, \ldots, e, f_0 \right) \alpha _0 = f_0,
\end{multline*}
because $\Remainderm{p - 1} = p - 1 < m - 1$ and $\Dividerm{p - 1} = 0$. Hence,
the relation $R_A \Pair{0}{p}$ is true. The following equality holds
\[
    f_1 ^{m} f_0 = \left( f_1, f_1, \ldots, f_1 \right) \sigma ^{0} \cdot
    \left( e, e, \ldots, e, f_0 \right) \alpha_0 = \left( f_1, \ldots, f_1, f_1
    f_0 \right) \alpha_0,
\]
whence
\[
    f_0 f_1 ^{m} f_0 = \left( e, \ldots, e, f_0 \right) \alpha _0 \cdot \left(
    f_1, \ldots, f_1, f_1 f_0 \right) \alpha _0 = f_1 ^m f_0,
\]
and $R_B \Argument{0}$ holds.

\par Now let $k \ge 1$. Using~\eqref{eq:unrolled_form_p_01}
and~\eqref{eq:unrolled_form_v_k}, decomposition of the left-hand part of the
relation of type A of length $k$ is defined by the following equality
\begin{multline*}
    f_0 f_1 ^{p m ^k - 1} v_k = \left( f_0 f_1 ^{p m ^{k - 1} - 1} v_{k - 1},
    f_0 f_1 ^{p m ^{k - 1} - 1} v_{k - 1}, \ldots, f_0 f_1 ^{p m ^{k - 1} - 1}
    v_{k - 1} \right) \alpha _0\\
    = \left( v_{k - 1}, v_{k - 1}, \ldots, v_{k - 1} \right) \alpha _0 = v_{k},
\end{multline*}
and the last equality is true due to the induction hypothesis
\[
    R_A \Pair{k - 1}{p} : \, f_0 f_1 ^{p m ^{k - 1} - 1} v_{k - 1} = v_{k - 1}.
\]
Hence, the relations $R_A \Pair{k}{p}$ hold in $\Semigroup{\Intermediate}$.
Similarly, let us write the decomposition of the left-hand part of the relation
$R_B \Argument{k}$:
\begin{multline*}
    f_0 f_1 ^{m ^k + m ^{k + 1} - 1} v_k = \\
    = \left( f_0 f_1 ^{ \left( m + 1 \right) m ^{k - 1} - 1} v_{k - 1}, f_0 f_1
    ^{ \left( m + 1 \right) m ^{k - 1} - 1} v_{k - 1}, \ldots, f_0 f_1 ^{
    \left( m + 1 \right) m ^{k - 1} - 1} v_{k - 1} \right) \alpha _0\\
    = \left( f_1 ^{m ^k} v_{k - 1}, f_1 ^{m ^k} v_{k - 1}, \ldots, f_1 ^{m
    ^k} v_{k - 1} \right) \alpha _0 = f_1 ^{m
    ^{k + 1}} v_{k},
\end{multline*}
where the equality of decompositions is substantiated by the induction
hypothesis for the relation $R_B \Argument{k - 1}$.
\end{proof}

\begin{proposition} \label{prop:all_relations}
In the semigroup $\Semigroup{\Intermediate}$ the relation
\begin{multline} \label{eq:all_relations}
    f_0 f_1 ^{m ^k p_{k + 2} - 1} \cdot f_1 ^{m ^{k + 1} p_{k + 1}} f_0 f_1
    ^{m ^{k} p_k - 1} f_0 f_1 ^{m ^{k - 1} p_{k - 1} - 1} f_0 \ldots f_1 ^{m
    p_1 - 1} f_0 \\
    = f_1 ^{m ^{k + 1} p_{k + 1}} f_0 f_1 ^{m ^{k} p_k - 1} f_0 f_1 ^{m ^{k -
    1} p_{k - 1} - 1} f_0 \ldots f_1 ^{m p_1 - 1} f_0,
\end{multline}
where $k \ge 0$, $1 \le p_{k + 2} \le {m - 1}$, $p_{k + 1} \ge 0$, $p_i \ge 1$,
$i = 1, 2, \ldots, {k}$, follows from the set of relations
\begin{align} \label{eq:defining_relations}
    R_A & \Pair{k}{p}, k \ge 0, p = 1, 2, \ldots, {m - 1}, & R_B &
    \Argument{k}, k \ge 0.
\end{align}
\end{proposition}

\begin{remark}
Let us denote the relation~\eqref{eq:all_relations} for fixed values of $k$,
$p_1, p_2, \ldots, p_{k + 2}$ by the symbol $r \left( k, p_{k + 2}, p_{k + 1},
p_{k}, \ldots, p_{1} \right)$, and we call $k$ as ``the length of this
relation''. In addition, the relations of type A and B can be written in the
form~\eqref{eq:all_relations}, because $R_A \Pair{k}{p} = r \left( k, p, 0, 1,
1, \ldots, 1 \right)$ and $R_B \Argument{k} = r \left( k, 1, 1, 1, \ldots, 1
\right)$.
\end{remark}

\begin{proof}
Let us prove the lemma by induction on $k$. For $k = 0$ the
relation~\eqref{eq:all_relations} is defined by the following equality
\begin{equation*}
    f_0 f_1 ^{p_{2} - 1} \cdot f_1 ^{m p_{1}} f_0 = f_1 ^{m p_{1}} f_0,
\end{equation*}
where $p_1 \ge 0$, $1 \le p_2 \le {m - 1}$. Using the relation $R_B
\Argument{0} : \, f_0 f_1 ^{m} f_0 = f_1 ^{m} f_0$, for any $p \ge 1$ the
following equalities hold
\begin{multline*}
    f_0 f_1 ^{mp} f_0 = f_0 f_1 ^{m \Argument{p - 1}} \cdot f_0 f_1 ^{m} f_0 =
    \ldots = f_0 f_1 ^m f_0 \left( f_1 ^m f_0 \right) ^{p - 1}\\
    = f_1 ^m \cdot f_0 f_1 ^m f_0 \left( f_1 ^m f_0 \right) ^{p - 2} = f_1
    ^{2m} f_0 \left( f_1 ^m f_0 \right) ^{p - 2} = \ldots = f_1 ^{mp} f_0.
\end{multline*}
Using the equality~$f_0 f_1 ^{m p_1} f_0 = f_1 ^{m p_1} f_0$ and the
relation~$R_A \Pair{0}{p_2}$ we have
\[
    f_0 f_1 ^{p_{2} + m p_{1} - 1} f_0 = f_0 f_1 ^{p_{2} - 1} f_0 f_1 ^{m
    p_{1}} f_0 = f_0 f_1 ^{m p_{1}} f_0 = f_1 ^{m p_{1}} f_0,
\]
whence the relation $r \left( 0, p_2, p_1 \right)$ holds, and is output from
the set~\eqref{eq:defining_relations}.

\par Let $k \ge 1$, and $p_1, p_2, \ldots, p_{k + 2}$ be integers that
fulfill the requirements of the lemma. Any relation~\eqref{eq:all_relations} of
length $\Argument{k - 1}$ is output from the set~\eqref{eq:defining_relations}
by induction hypothesis, and now we show that the
relation~\eqref{eq:all_relations} of length $k$ is output from the
relations~\eqref{eq:defining_relations} and the
relations~\eqref{eq:all_relations} of length $\Argument{k - 1}$.

\par Let $p_{k + 1} \ge 0$, $p_i \ge 1$, $i = 1, 2, \ldots, k$, be arbitrary
integers, and let us denote $w_p = f_1 ^{m ^{k + 1} p_{k + 1}} f_0 f_1 ^ {m
^{k} p_k - 1} f_0 f_1 ^ {m ^{k - 1} p_{k - 1} - 1} f_0 \ldots f_1 ^ {m p_1 - 1}
f_0$. Below we prove that the following equality holds
\begin{equation} \label{eq:insert_v_k}
    w_p = v_{k} \cdot w_p.
\end{equation}
Then the relation $r \left( k, p_{k + 2}, p_{k + 1}, \ldots, p_{1} \right)$
immediately follows from the equality \eqref{eq:insert_v_k} and~$R_A
\Pair{k}{p_{k + 2}}$:
\begin{equation*}
    f_0 f_1 ^{m ^k p_{k + 2} - 1} \cdot w_p = f_0 f_1 ^{m ^k p_{k + 2} - 1} v_k
    \cdot w_p = v_k \cdot w_p = w_p.
\end{equation*}

\par In order to prove~\eqref{eq:insert_v_k} we show that for any $p_{k}
\ge 0$, $p_i \ge 1$, $i = 1, 2, \ldots, {k - 1}$, the following equality holds
\begin{multline} \label{eq:insert_v_k-1}
    f_1 ^{m ^{k} p_{k}} f_0 f_1 ^ {m ^{k - 1} p_{k - 1} - 1} f_0 f_1 ^ {m ^{k -
    2} p_{k - 2} - 1} f_0 \ldots f_1 ^ {m p_1 - 1} f_0\\
    = v_{k - 1} \cdot f_1 ^{m ^{k} p_{k}} f_0 f_1 ^{m ^{k - 1} p_{k - 1} - 1}
    f_0 f_1 ^{m ^{k - 2} p_{k - 2} - 1} f_0 \ldots f_1 ^{m p_1 - 1} f_0,
\end{multline}
and then prove~\eqref{eq:insert_v_k} by induction on $p_{k + 1}$. We have
\begin{align*}
    f_1 ^{m ^{k} p_{k}} & f_0 f_1 ^{m ^{k - 1} p_{k - 1} - 1} f_0 f_1 ^{m ^{k -
    2} p_{k - 2} - 1} f_0 \ldots f_1 ^{m p_1 - 1} f_0 = \\
    & \begin{aligned}
    = f_0 f_1 ^{m ^{k - 1} - 1} \cdot f_1 ^{m ^{k - 1} \Argument{m p_{k}}} f_0
    f_1 ^ {m ^{k - 2} \Argument{m p_{k - 1}} - 1} f_0 & f_1 ^{m ^{k - 3}
    \Argument{m p_{k - 2}} - 1} f_0 \ldots\\
    & f_1 ^{m \Argument{m p_2} - 1} f_0 f_1 ^{m p_1 - 1} f_0
    \end{aligned}\\
    & \begin{aligned}
    = f_0 f_1 ^{m ^{k - 1} - 1} f_0 f_1 ^{m ^{k - 2} - 1} \cdot f_1 ^{m ^{k -
    2} \Argument{m ^2 p_{k}}} & f_0 f_1 ^ {m ^{k - 3} \Argument{m ^2 p_{k - 1}} -
    1} f_0 \ldots\\
    & f_1 ^{m \Argument{m ^2 p_3} - 1} f_0 f_1 ^{m ^2 p_2 - 1} f_0 f_1
    ^{m p_1 - 1} f_0
    \end{aligned}\\
    & = \ldots\\
    & \begin{aligned}
    = f_0 f_1 ^{m ^{k - 1} - 1} f_0 & f_1 ^ {m ^{k - 2} - 1} f_0 \ldots f_1 ^{m
    - 1} f_0 \cdot\\
    & \cdot f_1 ^{m ^{k} p_{k}} f_0 f_1 ^{m ^{k - 1} p_{k - 1} - 1} f_0 f_1 ^{m ^{k -
    2} p_{k - 2} - 1} f_0 \ldots f_1 ^{m p_1 - 1} f_0,
    \end{aligned}
\end{align*}
where each expansion of a semigroup word is substantiated by application of the
relation
\[
    r \left( k - i - 1, 1, m ^i p_{k}, m ^i p_{k - 1}, \ldots, m ^i p_{i + 1}
    \right)
\]
for $i = 0, 1, \ldots, {k - 1}$.

\par Now we prove~\eqref{eq:insert_v_k}, and let $p_{k + 1} = 0$.
Applying~\eqref{eq:insert_v_k-1}, ``reversed'' relation $R_A \Pair{k}{1} : v_k
= f_0 f_1 ^{m^k - 1} v_k$ and again the equality~\eqref{eq:insert_v_k-1}, the
following equalities hold
\begin{align*}
    f_0 & f_1 ^{m ^{k} p_k - 1} f_0 f_1 ^{m ^{k - 1} p_{k - 1} - 1} f_0 \ldots
    f_1 ^{m p_1 - 1} f_0 = \\
    & = f_0 f_1 ^{m ^{k} - 1} \cdot v_{k - 1} f_1 ^{m ^{k} \Argument{p_k - 1}}
    f_0 f_1 ^{m ^{k - 1} p_{k - 1} - 1} f_0 \ldots f_1 ^{m p_1 - 1} f_0 \\
    & = v_k \cdot f_1 ^{m ^{k} \Argument{p_k - 1}} f_0 f_1 ^{m ^{k - 1} p_{k -
    1} - 1} f_0 \ldots f_1 ^{m p_1 - 1} f_0 \\
    & = f_0 f_1 ^{m ^k - 1} v_k \cdot f_1 ^{m ^{k} \Argument{p_k - 1}} f_0 f_1
    ^{m ^{k - 1} p_{k - 1} - 1} f_0 \ldots f_1 ^{m p_1 - 1} f_0 \\
    & = f_0 f_1 ^{m ^k - 1} \cdot f_0 f_1 ^{m ^{k - 1} \Argument{m p_k} - 1}
    f_0 f_1 ^{m ^{k - 2} \Argument{m p_{k - 1}} - 1} f_0 \ldots f_1 ^{m
    \Argument{m p_2} - 1} f_0 f_1 ^{m p_1 - 1} f_0 \\
    & = f_0 f_1 ^{m ^k - 1} \cdot v_{k - 1} f_0 f_1 ^{m ^{k} p_k - 1} f_0 f_1
    ^{m ^{k - 1} p_{k - 1} - 1} f_0 \ldots f_1 ^{m p_1 - 1} f_0 \\
    & = v_{k} \cdot f_0 f_1 ^{m ^{k} p_k - 1} f_0 f_1 ^{m ^{k - 1} p_{k - 1} -
    1} f_0 \ldots f_1 ^{m p_1 - 1} f_0.
\end{align*}

\par Let $p_{k + 1} \ge 1$. The induction hypothesis is used for the adding
$v_k$, and the relation~$R_B \Argument{k}$ allows to add the word $f_0 f_1 ^{m
^k - 1}$. Then the word $v_k$ is cancelled, and the
equality~\eqref{eq:insert_v_k-1} is applied. Thus, the following equalities
hold
\begin{align*}
    f_1 ^{m ^{k + 1} p_{k + 1}} & f_0 f_1 ^{m ^{k} p_k - 1} f_0 f_1 ^{m ^{k -
    1} p_{k - 1} - 1} f_0 \ldots f_1 ^{m p_1 - 1} f_0 = \\
    & = f_1 ^{m ^{k + 1}} \cdot v_k \cdot f_1 ^{m ^{k + 1} \Argument{p_{k + 1}
    - 1}} f_0 f_1 ^{m ^{k} p_k - 1} f_0 f_1 ^{m ^{k - 1} p_{k - 1} - 1} f_0
    \ldots f_1 ^{m p_1 - 1} f_0\\
    & = f_0 f_1 ^{m ^k - 1} \cdot f_1 ^{m ^{k + 1}} v_k \cdot f_1 ^{m ^{k + 1}
    \Argument{p_{k + 1} - 1}} f_0 f_1 ^{m ^{k} p_k - 1} f_0 \ldots f_1 ^{m p_1
    - 1} f_0\\
    & = f_0 f_1 ^{m ^k - 1} \cdot f_1 ^{m ^{k + 1} p_{k + 1}} f_0 f_1 ^{m ^{k}
    p_k - 1} f_0 f_1 ^{m ^{k - 1} p_{k - 1} - 1} f_0 \ldots f_1 ^{m p_1 - 1}
    f_0 \\
    & = f_0 f_1 ^{m ^{k} - 1} \cdot v_{k - 1} \cdot f_1 ^{m ^{k + 1} p_{k + 1}}
    f_0 f_1 ^{m ^{k} p_k - 1} f_0 f_1 ^{m ^{k - 1} p_{k - 1} - 1} f_0 \ldots
    f_1 ^{m p_1 - 1} f_0 \\
    & = v_k \cdot f_1 ^{m ^{k + 1} p_{k + 1}} f_0 f_1 ^{m ^{k} p_k - 1} f_0 f_1
    ^{m ^{k - 1} p_{k - 1} - 1} f_0 \ldots f_1 ^{m p_1 - 1} f_0.
\end{align*}
The proposition is completely proved.
\end{proof}

\subsection{Reducing of semigroup words} \label{subsect:semigroup_reducing}

The main result of this subsection is the following proposition
\begin{proposition} \label{prop:normal_form}
Each element $\s \in \Semigroup{\Intermediate}$ can be reduced to the following
form
\begin{equation} \label{eq:normal_form}
    f_1 ^{p_k} f_0 f_1 ^{m ^{k - 1} p_{k - 1} - 1} f_0 \ldots f_1 ^{m ^{i}
    p_{i} - 1} f_0 \ldots f_1 ^{m ^2 p_{2} - 1} f_0 f_1 ^{m p_{1} - 1} f_0
    f_1 ^{p_0}
\end{equation}
where $k \ge 0$, $p_0 \ge 0$, $p_k \ge 0$, $p_i \ge 1$, $i = 1, 2, \ldots, {k -
1}$. There exists the reducing algorithm with complexity $\mathcal{O} \left(
\left| \s \right| \log_m \left| \s \right| \right)$.
\end{proposition}

\begin{remark}
In further we call the form~\eqref{eq:normal_form} as the (normal) form of
length $k$.
\end{remark}

\par It follows from Proposition~\ref{prop:all_relations}
that any relation~\eqref{eq:all_relations} cancels the beginning $f_0 f_1 ^p$
of a semigroup word for some $p$. Hence the reducing algorithm may run through
a semigroup word from the right-hand to the left-hand side, and it finishes
when reaches the beginning of $\s$ (or the most right symbol $f_0$). In this
subsection we consider the reducing of a semigroup word written in special
form, and then describe the reducing algorithm. The proof of
Proposition~\ref{prop:normal_form} bases on these results.

\par Let $\s$ be an arbitrary semigroup word such that
\begin{equation*}
    \s = f_0 f_1 ^{m ^{k - 1} p_{k - 1} - 1} f_0 f_1 ^{m
    ^{k - 2} p_{k - 2} - 1} \ldots f_0 f_1 ^{m p_1 - 1} f_0 f_1 ^{p_0},
\end{equation*}
where $k \ge 1$, $p_{i} \ge 1$, $i = 1, 2, \ldots, k - 1$, $p_0 \ge 0$, and let
us consider the following semigroup word
\begin{equation*}
    \s' = f_0 f_1 ^{p_k} \s = f_0 f_1 ^{p_k} f_0 f_1 ^{m ^{k - 1} p_{k - 1} -
    1} f_0 f_1 ^{m ^{k - 2} p_{k - 2} - 1} f_0 \ldots f_1 ^{m p_1 - 1} f_0 f_1
    ^{p_0},
\end{equation*}
where $p_k \ge 1$. It follows from Proposition~\ref{prop:all_relations} that
the relations~\eqref{eq:all_relations} can be applied to $\s'$, if there exist
the integers $0 \le i \le {k - 1}$, and $q_0 \ge 0$, $q_1 \in \left\{ 1, 2,
\ldots, {m - 1} \right\}$ such that $p_k$ can be presented by the equality
\begin{equation} \label{eq:reduced_p_k}
    p_k = m ^{i} q_1 + m ^{i + 1} q_0 - 1.
\end{equation}
Then the relation
\[
    r \left( i, q_1, q_0, m ^{k - 1 - i} p_{k - 1}, m ^{k - 1 - i} p_{k - 2},
    \ldots, m ^{k - 1 - i} p_{k - i} \right)
\]
can be used in order to cancel the beginning $f_0 f_1 ^{m ^{i} q_1 - 1}$ of
$\s'$. Clearly $q_0$ and $q_1$ are unambiguously defined by $p_k$.

\par Let $p \ge 1$ be an arbitrary integer, and let us denote
\[
    t_1 \n[p] = \max \left\{ {\begin{array}{*{20}c}
        {j \ge 0} & \vline & {m ^j \; \vline \; p } \\
    \end{array}} \right\},
\]
that is the maximal power of $m$ such that $p$ is divisible by $m ^{t_1
\n[p]}$. Similarly, let $t_2 \n[p]$ is defined by the equality
\[
    t_2 \n[p] = p \mod m ^{t_1 \n[p] + 1}.
\]
Obviously for any $p \ge 1$ the number ${t_2 \n[p]} / {m ^{t_1 \n[p]}}$ is the
positive integer such that
\[
    1 \le \frac{t_2 \n[p]}{m ^{t_1 \n[p]}} < m.
\]
Using these definitions, the integer $p$ can be written as
\[
    p = m ^{t_1 \n[p]} \frac{p}{m ^{t_1 \n[p]}} = m ^{t_1 \n[p] + 1} \left[
    \frac{p}{m ^{t_1 \n[p] + 1}} \right] + m ^{t_1 \n[p]} \frac{t_2 \n[p]}{m
    ^{t_1 \n[p]}}.
\]

\par If we assume
\begin{align*}
    q_0 & = \left[ \frac{p_k + 1}{m ^{t_1 \n[p_k + 1] + 1}} \right] & &
    \text{and} & q_1 & = \frac{t_2 \n[p_k + 1]}{m ^{t_1 \n[p_k + 1]}},
\end{align*}
then $q_0 \ge 0$ and $1 \le q_1 \le {m - 1}$, and these numbers satisfy the
equality~\eqref{eq:reduced_p_k} for $i = t_1 \n[p_k + 1]$. If $t_1 \n[p_k + 1]
< {k}$, then the relation
\begin{multline*}
    r \left( t_1 \n[p_k + 1], \frac{t_2 \n[p_k + 1]}{m ^{t_1 \n[p_k + 1]}},
    \left[ \frac{p_k + 1}{m ^{t_1 \n[p_k + 1] + 1}} \right], m ^{k - 1 - t_1
    \n[p_k + 1]} p_{k - 1}, \right. \\
    \left. m ^{k - 1 - t_1 \n[p_k + 1]} p_{k - 2}, \ldots, m ^{k - 1 - t_1
    \n[p_k + 1]} p_{k - t_1 \n[p_k + 1]} \right)
\end{multline*}
allows to cancel the semigroup word
\[
    f_0 f_1 ^{t_2 \n[p_k + 1] - 1}
\]
at the beginning of $\s'$. Hence, the element $\s$ is equivalent to the
following element
\[
    \s = f_1 ^{m ^{t_1 \n[p_k + 1] + 1} \left[ \frac{p_k + 1}{m ^{t_1 \n[p_k +
    1] + 1}} \right]} f_0 f_1 ^{m ^{k - 1} p_{k - 1} - 1} f_0 f_1 ^{m ^{k - 2}
    p_{k - 2} - 1} \ldots f_0 f_1 ^{m p_1 - 1} f_0 f_1 ^{p_0}.
\]

\incmargin{1em} \linesnumbered
\begin{algorithm}[t] \label{alg:reducing}
    \SetKwData{Word}{$\s$}
    \SetKwData{Index}{i}
    \SetKwData{RIndex}{j}
    \SetKwData{Degree}{r}

    \KwData{A semigroup word
    \begin{equation*}
        \Word = f_1 ^{p_{2k}} f_0 ^{p_{2k - 1}} f_1 ^{p_{2k - 2}} \ldots f_0
        ^{p_{1}} f_1 ^{p_{0}},
    \end{equation*}
    where $k \ge 0$, $p_0, p_{2k} \ge 0$, $p_i \ge 1$, $i = 1, 2, \ldots, {2k -
    1}$.}
    \KwResult{A semigroup word \Word written in the form~\eqref{eq:normal_form}.}
    \caption{The Reducing Algorithm}

    \BlankLine
    \Index $\longleftarrow 0$  \; \label{alg_reducing:init_begin}
    \RIndex $\longleftarrow 0$ \;
    \Degree $\longleftarrow 0$ \; \label{alg_reducing:init_end}

    \BlankLine
    \For{\Index $\longleftarrow 1$ \KwTo $2k - 1$ \label{alg_reducing:for_begin}}{
        \eIf{\Index is odd}
        {
            $p_i \longleftarrow 1$ \;
            \label{alg_reducing:if_odd_begin}
            $\RIndex \longleftarrow \RIndex + 1$ \;
            $r \longleftarrow 0$ \;
            \label{alg_reducing:if_odd_end}
        } 
        {
            \eIf{$\left( p_i + r \right) \mod m ^j = m ^j - 1$ \label{alg_reducing:if_even_begin}}
            {
                $p_i \longleftarrow \Argument{p_i + r}$ \; \label{alg_reducing:if_even1}
            }
            {
                The subword $f_0 ^{p_{i + 1}} f_1 ^{p_i}$ is canceled in \Word
                \;
                \label{alg_reducing:if_even2_begin}
                $\Degree \longleftarrow p_i + \Degree - t_2 \n[p_i + \Degree +
                1] + 1$ \;
                \label{alg_reducing:if_even2_end}
                $\Index \longleftarrow \Index + 1$ \;
                \label{alg_reducing:if_even_end}
            }
        } 
     } \label{alg_reducing:for_end}
\end{algorithm}

\begin{proof}[Proof of Proposition~\ref{prop:normal_form}]
Let us consider Algorithm~\ref{alg:reducing}. We prove that it reduces an
arbitrary semigroup word $\s$ to the form~\eqref{eq:normal_form}.

\par The local variables are initialized at
lines~\ref{alg_reducing:init_begin}--\ref{alg_reducing:init_end}, and it is
executed once. There $\mathsf{i}$ is the index of exponent in the input word
$\s$, $\mathsf{j}$ is the index of exponent in reduced part of the semigroup
word, and $\mathsf{r}$ is a temporary variable, that is used for calculating
the values of exponents in the reduced word.

\par The main loop at
lines~\ref{alg_reducing:for_begin}--\ref{alg_reducing:for_end} moves along $\s$
from the right-hand side, and sequentially reduces exponents at the symbol
$f_1$ to the form $m ^j q_j - 1$, where $q_j > 0$ and $j$ varies over the
values $0, 1, 2, \ldots$. If $k = 0$ and $\s = f_1 ^{p_0}$, then $\s$ is
already of the form~\eqref{eq:normal_form}. In this case the main loop is not
executed. Otherwise, let us consider the $i$-th iteration of the main loop,
where the algorithm checks the value of $p_i$.

\par If $i$ is odd, then the
lines~\ref{alg_reducing:if_odd_begin}--\ref{alg_reducing:if_odd_end} is
executed. In this case $p_i$ is the exponent at the symbol $f_0$, and the
subword $f_0 ^{p_i - 1}$ can be canceled by the applying the relation $f_0 ^2 =
f_0$. Therefore $p_i$ is assigned to $1$, the algorithm starts ``to collect''
the next exponent at $f_1$ in the reduced word, and the loop moves to the next
value of $i$.

\par If $i$ is even, then $p_i$ is exponent at $f_1$. At the
line~\ref{alg_reducing:if_even_begin} the semigroup word $\s$ is defined by the
following equality
\[
    f_1 ^{p_{2k}} f_0 ^{p_{2k - 1}} f_1 ^{p_{2k - 2}} \ldots f_0
    ^{p_{\mathsf{i} + 1}} \cdot f_1 ^{p_{\mathsf{i}}} f_1 ^{r} \cdot f_0 f_1
    ^{m ^{j - 1} p' _{j - 1} - 1} \ldots f_0 f_1 ^{m p' _1 - 1} f_0 f_1 ^{p_0},
\]
where $p' _q \ge 1$, $q = 1, 2, \ldots, {j - 1}$, and $p_{\mathsf{i}} \ge 1$,
$r \ge 0$; and let us separate it into two parts
\begin{align*}
    \s_1 & = f_1 ^{p_{2k}} f_0 ^{p_{2k - 1}} f_1 ^{p_{2k - 2}} \ldots f_0
    ^{p_{{\mathsf{i}} + 1} - 1},\\
    \s_2 & = f_0 f_1 ^{p_{\mathsf{i}} + r} \cdot f_0 f_1 ^{m ^{j - 1} p' _{j -
    1} - 1} \ldots f_0 f_1 ^{m p' _1 - 1} f_0 f_1 ^{p_0}.
\end{align*}
If the equality $p_{\mathsf{i}} + r = m ^j p' _j - 1$ holds for some $p'_j >
0$, then $\s_2$ has already written in the form~\eqref{eq:normal_form}. Then
the line~\ref{alg_reducing:if_even1} is executed, and the algorithm continues
on the next exponent of $\s$.

\par Otherwise, it follows from the speculations above that
$\s_2$ is reducible. The subword $f_0 f_1 ^{t_2 \n[p_{\mathsf{i}} + r + 1] -
1}$ is cancelled, and the subword $f_0 ^{p_{{\mathsf{i}} + 1} - 1}$ is
cancelled due to the relation $f_0 ^2 = f_0$. Therefore the algorithm cancels
the subword $f_0 ^{p_{\mathsf{i} + 1}} f_1 ^{p_{\mathsf{i}}}$ at the
line~\ref{alg_reducing:if_even2_begin}, but increases $r$ at the next line.
Then the loop continues on the exponent $p_{\mathsf{i} + 2}$ at the next symbol
$f_1$.

\par The number of iterations of the main loop is equal to $2k - 1$, where $k$ is
defined by the input word. Clearly $2k \le \left| \s \right|$. Each iteration
includes fixed number of arithmetic and logical operations, and calculating of
$t_2$. As $\left( p_{\mathsf{i}} + r + 1 \right) < \left| \s \right|$, thus the
complexity of $t_2 \left( p_{\mathsf{i}} + r + 1 \right)$ calculating is not
greater than $\log _m \left| \s \right|$. Therefore there exists the positive
integer $c_1$ such that the complexity of one main loop iteration doesn't
exceed $c_1 + \log _m \left| \s \right|$, whence the total complexity of
Algorithm~\ref{alg:reducing} equals $\mathcal{O} \Argument{\left| \s \right|
\log _m \left| \s \right|}$. Obviously the real complexity depends on algorithm
realization.
\end{proof}

\subsection{Normal form} \label{subsect:semigroup_normal_form}

\par It follows from the previous subsection that each element can be reduced to
the form~\eqref{eq:normal_form}. The main result of this subsection is that two
semigroup elements written in different form~\eqref{eq:normal_form} define
different automaton transformations. Namely,

\begin{proposition} \label{prop:normality_of_form}
\par Let  $\s _1$, $\s _2$ be arbitrary elements of the semigroup
$\Semigroup{\Intermediate}$ written in the form~\eqref{eq:normal_form}:
\begin{equation*}
\begin{aligned}
    \s_1 & = f_1 ^{p_k} f_0 f_1 ^{m ^{k - 1} p_{k - 1} - 1} f_0 \ldots f_1 ^{m
    ^2 p_{2} - 1} f_0 f_1 ^{m p_{1} - 1} f_0 f_1 ^{p_0},\\
    \s_2 & = f_1 ^{q_l} f_0 f_1 ^{m ^{l - 1} q_{l - 1} - 1} f_0 \ldots f_1 ^{m
    ^2 q_{2} - 1} f_0 f_1 ^{m q_{1} - 1} f_0 f_1 ^{q_0},
\end{aligned}
\end{equation*}
where $k \ge 0$, $l \ge 0$, $p_0, p_k \ge 0$, $q_0, q_l \ge 0$, $p_i \ge 1$, $i
= 1, 2, \ldots, {k - 1}$, $q_j \ge 1$, $j = 1, 2, \ldots, {l - 1}$. Then $\s_1$
and $\s_2$ define the same automaton transformation over $\InfiniteSet$ if and
only if they coincide graphically, that is
\begin{equation*}
    k = l, p_0 = q_0, p_1 = q_1, \ldots, p_k = q_l.
\end{equation*}
\end{proposition}

\par Before the proof we consider the restrictions of arbitrary semigroup
element written in the form~\eqref{eq:normal_form}. Let us introduce two
functions $r_1, r_2 : \Natural \to \left\{ 0, 1, \ldots, m - 1 \right\}$ such
that for any $p \in \Natural$ they are defined by the equalities
\begin{align*}
    r_1 \n[p] & = \delta _{m - 1, \Remainderm{p}} = \left\{
        \begin{array}{ll}
            0, & \hbox{if $0 \le \Remainderm{p} \le {m - 2}$,} \\
            1, & \hbox{if $\Remainderm{p} = {m - 1}$;} \\
        \end{array}
    \right.\\
    r_2 \n[p] & = m - 1 - \Remainderm{p}.
\end{align*}
Clearly for any $p \in \Natural$ the inequality $r_1 \n[p] \neq r_2 \n[p]$
holds.

\par Let $\s \in \Semigroup{\Intermediate}$ be a semigroup element written in
the form~\eqref{eq:normal_form} of length $k = 1$:
\[
    \s = f_1 ^{p_1} f_0 f_1 ^{p_0},
\]
where $p_0, p_1 \ge 0$. It follows from Lemma~\ref{lem:power_f_1}
and~\eqref{eq:unrolled_form_p_01} that $\s$ has the following decomposition
\[
    \s = \left( f_1 ^{p}, \ldots, f_1 ^{p}, f_1 ^{\Dividerm{p_1}} f_0 f_1
    ^{\Dividerm{p_0}}, f_1 ^{p + 1}, \ldots, f_1 ^{p + 1} \right) \alpha
    _{\Remainder{p_1}},
\]
where $p = \Dividerm{p_1} + \Dividerm{p_0}$, whence
\begin{subequations}
\begin{align}
    \label{eq:restriction_k=1_r_1}
    \Restriction{\s}{x_{r_1 \n[p_0]}} & = f_1 ^{\Dividerm{p_1} + \Dividerm{p_0} +
    r_1 \n[p_0]},\\
    \label{eq:restriction_k=1_r_2}
    \Restriction{\s}{x_{r_2 \n[p_0]}} & = f_1 ^{\Dividerm{p_1}} f_0 f_1
    ^{\Dividerm{p_0}},\\
\intertext{and for all $0 \le r \le {m - 1}$, $r \neq r_2 \n[p_0]$, the
equality hold}
    \label{eq:restriction_k=1_r}
    \Restriction{\s}{x_{r}} & = f_1 ^{\Dividerm{p_1} + \Dividerm{p_0 + r}}.
\end{align}
All elements $\Restriction{\s}{x_{r_1 \n[p_0]}}$, $\Restriction{\s}{x_{r_2
\n[p_0]}}$, and $\Restriction{\s}{x_{r}}$ are written in the
form~\eqref{eq:normal_form}.
\end{subequations}

\par Now let $\s \in \Semigroup{\Intermediate}$ be a semigroup
element written in the form~\eqref{eq:normal_form}:
\[
    \s = f_1 ^{p_k} f_0 f_1 ^{m ^{k - 1} p_{k - 1} - 1} f_0 \ldots f_1 ^{m ^{i}
    p_{i} - 1} f_0 \ldots f_1 ^{m ^2 p_2 - 1} f_0 f_1 ^{m p_1 - 1} f_0 f_1 ^{p_0},
\]
where $k > 1$, $p_0 \ge 0$, $p_k \ge 0$, $p_i \ge 1$, $i = 1, 2, \ldots, {k -
1}$. It follows from Lemma~\ref{lem:power_f_1}
and~\eqref{eq:unrolled_form_p_01} that $\s$ has the following decomposition
\begin{equation*}
\begin{aligned}
    \s = & \left( f_1 ^{\Dividerm{p_k}} f_0 f_1 ^{m^{k - 2} p_{k - 1} - 1}
    \ldots f_0 f_1 ^{m p_2 - 1} f_0 f_1 ^{p_{1} - 1} \cdot f_1
    ^{\Dividerm{p_0}}, \right.\\
    &  \ldots\\
    & f_1 ^{\Dividerm{p_k}} f_0 f_1 ^{m^{k - 2} p_{k - 1} - 1} \ldots f_0 f_1 ^{m
    p_2 - 1} f_0 f_1 ^{p_{1} - 1} \cdot f_1 ^{\Dividerm{p_0 + \left( m - 2 -
    \Remainderm{p_0} \right)}},\\
    & f_1 ^{\Dividerm{p_k}} f_0 f_1 ^{m^{k - 2} p_{k - 1} - 1}
    \ldots f_0 f_1 ^{m p_2 - 1} f_0 f_1 ^{p_{1} - 1} \cdot f_0 f_1
    ^{\Dividerm{p_0 + \left( m - 1 - \Remainderm{p_0} \right)}},\\
    & f_1 ^{\Dividerm{p_k}} f_0 f_1 ^{m^{k - 2} p_{k - 1} - 1}
    \ldots f_0 f_1 ^{m p_2 - 1} f_0 f_1 ^{p_{1} - 1} \cdot f_1
    ^{\Dividerm{p_0 + \left( m - \Remainderm{p_0} \right)}},\\
    & \ldots,\\
    & \left. f_1 ^{\Dividerm{p_k}} f_0 f_1 ^{m ^{k - 2} p_{k - 1} - 1} \ldots
    f_0 f_1 ^{m p_{2} - 1} f_0 f_1 ^{p_{1} - 1} \cdot f_1 ^{\Dividerm{p_0 + m -
    1}} \right) \alpha _{\Remainder{p_k}}.
\end{aligned}
\end{equation*}
Hence, the restrictions of $\s$ are defined by the following equalities
\begin{subequations}
\begin{align}
    \label{eq:restriction_k>1_r_1}
    \Restriction{\s}{x_{r_1 \n[p_0]}} & = f_1 ^{\Dividerm{p_k}} f_0 f_1 ^{m ^{k
    - 2} p_{k - 1} - 1} \ldots f_0 f_1 ^{m p_{2} - 1} f_0 f_1 ^{p_{1} - 1 +
    \Dividerm{p_{0}} + r_1 \n[p_0]},\\
    \label{eq:restriction_k>1_r_2}
    \Restriction{\s}{x_{r_2 \n[p_0]}} & = f_1 ^{\Dividerm{p_k}} f_0 f_1 ^{m ^{k
    - 2} p_{k - 1} - 1} \ldots f_0 f_1 ^{m p_{2} - 1} f_0 f_1 ^{p_{1} - 1} f_0
    f_1 ^{\Dividerm{p_0}},
\intertext{and}
    \label{eq:restriction_k>1_r}
    \Restriction{\s}{x_{r}} & = f_1 ^{\Dividerm{p_k}} f_0 f_1 ^{m ^{k - 2} p_{k -
    1} - 1} \ldots f_0 f_1 ^{m p_{2} - 1} f_0 f_1 ^{p_{1} - 1 + \Dividerm{p_0 +
    r}},
\end{align}
for any $0 \le r \le {m - 1}$, $r \neq r_2 \n[p_0]$. The elements
$\Restriction{\s}{x_{r_1 \n[p_0]}}$ and $\Restriction{\s}{x_{r}}$ are already
written in the form~\eqref{eq:normal_form} and are irreducible. On the other
hand, the semigroup word $\Restriction{\s}{x_{r_2 \n[p_0]}}$ may be reduced. If
all integers $p_1, p_2, \ldots, p_{k - 1}$ are divisible by $m$, then
$\Restriction{\s}{x_{r_2 \n[p_0]}}$ can be written in the
form~\eqref{eq:normal_form}:
\[
    \Restriction{\s}{x_{r_2 \n[p_0]}} = f_1 ^{\Dividerm{p_k}} f_0 f_1 ^{m ^{k -
    1} \Dividerm{p_{k - 1}} - 1} \ldots f_0 f_1 ^{m^2 \Dividerm{p_{2}} - 1} f_0 f_1
    ^{m \Dividerm{p_{1}} - 1} f_0 f_1 ^{\Dividerm{p_0}}.
\]
Otherwise, let $i_0$, $1 \le i_0 \le {k - 1}$, be the minimal index such that
$p_{i_0}$ is not divisible by $m$. Then the element $\Restriction{\s}{x_{r_2
\n[p_0]}}$ is reduced to the following element
\begin{multline*}
    \Restriction{\s}{x_{r_2 \n[p_0]}} = f_1 ^{\Dividerm{p_k}} f_0 f_1 ^{m^{k -
    2} p_{k - 1} - 1} \ldots f_0 f_1 ^{m ^{i_0 + 1} p_{i_0 + 2} - 1} f_0 f_1
    ^{m ^{i_0} \left( p_{i_0 + 1} + \Dividerm{p_{i_0}} \right) - 1} \\
    \cdot f_0 f_1 ^{m ^{i_0 - 1} \Argument{\frac{p_{i_0 - 1}}{m}} - 1} \ldots
    f_0 f_1 ^{m^2 \Argument{\frac{p_{2}}{m}} - 1} f_0 f_1 ^{m
    \Argument{\frac{p_{1}}{m}} - 1} f_0 f_1 ^{\Dividerm{p_0}}.
\end{multline*}
\end{subequations}

\begin{proof}[Proof of Proposition~\ref{prop:normality_of_form}]
Not restricting generality, let $0 \le k \le l$. Let us assume that the
elements $\s_1$ and $\s_2$ define the same automaton transformation over
$\InfiniteSet$. Then for any $u \in \InfiniteSet$ the equality holds
\begin{equation} \label{eq:assumption}
    \s_1 \Argument{u} = \s_2 \Argument{u},
\end{equation}
whence for any $v \in \FiniteSet$ the restrictions of $\s_1$ and $\s_2$
coincide, i.e. for arbitrary $u \in \InfiniteSet$ the equality holds
\[
    \Restriction{\s_1}{v} \n[u] = \Restriction{\s_2}{v} \n[u].
\]
We prove the proposition by induction on $k$.

\par Let $k = 0$, and $\s_1 = f_1 ^{p_0}$. If $l > 0$ then the transformation
$\s_2$ includes $f_0$ and is not bijective. In the case $l = 0$ for input word
$u_0 = \eta \n[0] = x_0 ^{\ast}$ we have
\begin{align*}
    \s_1 \n[u_0] & = f_1 ^{p_0} \n[u_0] = \eta \n[p_0],\\
    \s_2 \n[u_0] & = f_1 ^{q_0} \n[u_0] = \eta \n[q_0].
\end{align*}
It follows from the assumption~\eqref{eq:assumption} that $\eta \n[p_0] = \eta
\n[q_0]$, and, consequently, $p_0 = q_0$. Thus for $k = 0$ it follows
from~\eqref{eq:assumption} that the requirements $l = 0$ and $p_0 = q_0$ should
be fulfilled.

\par Now let $k \ge 1$, and there are two possible cases: $\Remainderm{p_0}
\neq \Remainderm{q_0}$ and $\Remainderm{p_0} = \Remainderm{q_0}$.

\par 1. Let $\Remainderm{p_0} \neq \Remainderm{q_0}$. It follows
from~\eqref{eq:restriction_k=1_r}, \eqref{eq:restriction_k>1_r},
\eqref{eq:restriction_k=1_r_2}, and~\eqref{eq:restriction_k>1_r_2} that for the
input word $x_{r_2 \n[q_0]} u$, $u \in \InfiniteSet[m]$, the following
equalities hold
\begin{align*}
    \s_1 \n[x_{r_2 \n[q_0]} u] & = x_{\Remainderm{p_k}} \cdot
    \Restriction{\s_1}{x_{r_2 \n[q_0]}} \n[u],\\
    \s_2 \n[x_{r_2 \n[q_0]} u] & = x_{\Remainderm{q_l}} \cdot
    \Restriction{\s_2}{x_{r_2 \n[q_0]}} \n[u],
\end{align*}
where
\[
    \Restriction{\s_1}{x_{r_2 \n[q_0]}} = \left\{
        \begin{array}{ll}
            {f_1 ^{\Dividerm{p_1} + \Dividerm{p_0 + r_2 \n[q_0]}}}, & \hbox{if
            $k = 1$;} \\
            {\begin{aligned}
            f_1 ^{\Dividerm{p_k}} f_0 f_1 ^{m ^{k - 2} p_{k - 1} - 1} & \ldots \\
            \ldots f_0 & f_1 ^{m p_{2} - 1} f_0 f_1 ^{p_{1} - 1 + \Dividerm{p_{0} + r_2
            \n[q_0]}}
            \end{aligned}}, & \hbox{otherwise;}
            \\
        \end{array}%
    \right.
\]
and
\[
    \Restriction{\s_2}{x_{r_2 \n[q_0]}} = \left\{
        \begin{array}{ll}
            {f_1 ^{\Dividerm{q_1}} f_0 f_1 ^{\Dividerm{q_0}}}, & \hbox{if
            $l = 1$;} \\
            {f_1 ^{\Dividerm{q_l}} f_0 f_1 ^{m ^{l - 2} q_{l - 1} - 1} \ldots
            f_0 f_1 ^{m q_{2} - 1} f_0 f_1 ^{q_{1} - 1} f_0 f_1
            ^{\Dividerm{q_{0}}}}, & \hbox{otherwise.}
            \\
        \end{array}%
    \right.
\]
The element $\Restriction{\s_1}{x_{r_2 \n[q_0]}}$ is irreducible, and has the
normal form of length $\Argument{k - 1}$. By induction hypothesis the element
$\Restriction{\s_2}{x_{r_2 \n[q_0]}}$ should have the normal form of length
$\left( k - 1 \right)$, but $\Restriction{\s_2}{x_{r_2 \n[q_0]}}$ has the
form~\eqref{eq:normal_form} of length $l$ or $\Argument{l - 1}$. It follows
from the condition $l \ge k$ that $l = k$ and $\Restriction{\s_2}{x_{r_2
\n[q_0]}}$ is reducible.

\par In the case $l = 1$ the element $\Restriction{\s_2}{x_{r_2 \n[q_0]}}$ is
irreducible and has the normal form of length $1$ ($ > 0$), so $l > 1$ and
there exists the minimal index $j_0$, $1 \le j_0 \le k - 1$, such that
$q_{j_0}$ is not divisible by $m$. The element $\Restriction{\s_2}{x_{r_2
\n[q_0]}}$ is written in the following form
\begin{multline*}
    \Restriction{\s_2}{x_{r_2 \n[q_0]}} = f_1 ^{\Dividerm{q_k}} f_0 f_1 ^{m ^{k
    - 2} q_{k - 1} - 1} \ldots f_0 f_1 ^{m ^{j_0 + 1} q_{j_0 + 2} - 1} f_0 f_1
    ^{m ^{j_0} \left( q_{j_0 + 1} + \Dividerm{q_{j_0}} \right) - 1} \\
    \cdot f_0 f_1 ^{m ^{j_0 - 1} \left( \frac{q_{j_0 - 1}}{m} \right) - 1}
    \ldots f_0 f_1 ^{m ^2 \left( \frac{q_2}{m} \right) - 1} f_0 f_1 ^{m \left(
    \frac{q_1}{m} \right) - 1} f_0 f_1 ^{\Dividerm{q_0}}.
\end{multline*}
It follows from the assumption~\eqref{eq:assumption} that the following set of
requirements should be fulfilled
\begin{gather*}
    k = l > 1, \Remainderm{p_k} = \Remainderm{q_k}, \Dividerm{p_k} =
    \Dividerm{q_k}, p_{k - 1} = q_{k - 1}, \ldots, p_{j_0 + 2} = q_{j_0 + 2},
    \\
    p_{j_0 + 1} = q_{j_0 + 1} + \Dividerm{q_{j_0}}, p_{j_0} = \frac{q_{j_0 -
    1}}{m}, \ldots, p_{2} = \frac{q_{1}}{m}, \\
    p_1 - 1 + \Dividerm{p_0 + r_2 \n[q_0]} = \Dividerm{q_0}.
\end{gather*}

As the equality
\begin{multline*}
    \Dividerm{p_0 + r_2 \n[q_0]} = \Dividerm{m \Dividerm{p_0} +
    \Remainderm{p_0} + m - 1 - \Remainderm{q_0}} \\
    = \Dividerm{p_0} + \left\{%
        \begin{array}{ll}
            0, & \hbox{$\Remainderm{p_0} <
            \Remainderm{q_0}$} \\
            1, & \hbox{$\Remainderm{p_0} >
            \Remainderm{q_0}$} \\
        \end{array}%
    \right.
\end{multline*}
holds, then the set of requirements can be written in the following way
\begin{equation} \label{eq:conditions_j0}
\begin{gathered}
    k = l, p_k = q_k, p_{k - 1} = q_{k - 1}, \ldots, p_{j_0 + 2} = q_{j_0 + 2},
    \\
    p_{j_0 + 1} = q_{j_0 + 1} + \Dividerm{q_{j_0}}, p_{j_0} = \frac{q_{j_0
    - 1}}{m}, \ldots, p_{2} = \frac{q_{1}}{m},\\
    p_1 = \Dividerm{q_0} - \Dividerm{p_0} + \left\{%
        \begin{array}{ll}
            1, & \hbox{$\Remainderm{p_0} < \Remainderm{q_0}$;} \\
            0, & \hbox{$\Remainderm{p_0} > \Remainderm{q_0}$.} \\
        \end{array}%
    \right.
\end{gathered}
\end{equation}

\par Similar reasoning can be carried out for the input word $x_{r_2 \n[p_0]}$,
where the elements $\s_1$ and $\s_2$ are rearranged. Hence, there exists the
minimal index $i_0$, $1 \le i_0 \le k - 1$, such that $p_{i_0}$ is not
divisible by $m$, and the following set of requirements should be fulfilled
\begin{equation} \label{eq:conditions_i0}
\begin{gathered}
    k = l, p_k = q_k, p_{k - 1} = q_{k - 1}, \ldots, p_{i_0 + 2} = q_{i_0 + 2},
    \\
    p_{i_0 + 1} + \Dividerm{p_{i_0}} = q_{i_0 + 1}, \frac{p_{i_0 - 1}}{m} =
    q_{i_0}, \ldots, \frac{p_{1}}{m} = q_{2}, \\
    q_1 = \Dividerm{p_0} - \Dividerm{q_0} + \left\{%
        \begin{array}{ll}
            0, & \hbox{$\Remainderm{p_0} < \Remainderm{q_0}$;} \\
            1, & \hbox{$\Remainderm{p_0} > \Remainderm{q_0}$.} \\
        \end{array}%
    \right.
\end{gathered}
\end{equation}
Summarizing two last requirements of~\eqref{eq:conditions_j0}
and~\eqref{eq:conditions_i0} we have the following equality:
\[
    p_1 + q_1 = 1.
\]
This equality contradicts the requirements $k = l > 1$ and $p_1, q_1 \ge 1$.
Hence, the contradiction with the assumption~\eqref{eq:assumption} is obtained.

\par 2. Let $\Remainderm{p_0} = \Remainderm{q_0}$. It follows
from~\eqref{eq:restriction_k=1_r_1} and~\eqref{eq:restriction_k>1_r_1} that for
an arbitrary word $u \in \InfiniteSet$ the equalities hold
\begin{align*}
    \s_1 \n[x_{r_1 \n[p_0]} u] & = x_{\Remainderm{p_k}} \cdot
    \Restriction{\s_1}{x_{r_1 \n[p_0]}} \n[u],\\
    \s_2 \n[x_{r_1 \n[p_0]} u] & = x_{\Remainderm{q_l}} \cdot
    \Restriction{\s_2}{x_{r_1 \n[p_0]}} \n[u],
\end{align*}
where
\[
    \Restriction{\s_1}{x_{r_1 \n[p_0]}} = \left\{
        \begin{array}{ll}
            {f_1 ^{\Dividerm{p_1} + \Dividerm{p_0} + r_1 \n[p_0]}}, & \hbox{if
            $k = 1$;} \\
            {\begin{aligned}
            f_1 ^{\Dividerm{p_k}} f_0 f_1 ^{m ^{k - 2} p_{k - 1} - 1} & \ldots
            \\
            \ldots f_0 & f_1 ^{m p_{2} - 1} f_0 f_1 ^{p_{1} - 1 +
            \Dividerm{p_{0}} + r_1 \n[p_0]}
            \end{aligned}}, & \hbox{otherwise;}
            \\
        \end{array}%
    \right.
\]
and
\[
    \Restriction{\s_2}{x_{r_1 \n[p_0]}} = \left\{
        \begin{array}{ll}
            {f_1 ^{\Dividerm{q_1} + \Dividerm{q_0} + r_1 \n[p_0]}}, & \hbox{if
            $l = 1$;} \\
            {f_1 ^{\Dividerm{q_l}} f_0 f_1 ^{m ^{l - 2} q_{l - 1} - 1} \ldots
            f_0 f_1 ^{m q_{2} - 1} f_0 f_1 ^{q_{1} - 1 + \Dividerm{q_{0}} + r_1
            \n[p_0]}}, & \hbox{otherwise.}
            \\
        \end{array}%
    \right.
\]
As $\Restriction{\s_1}{x_{r_1 \n[p_0]}}$ and $\Restriction{\s_2}{x_{r_1
\n[p_0]}}$ are written in the form~\eqref{eq:normal_form} and their normal form
has length of $\Argument{k - 1}$ and $\Argument{l - 1}$, respectively, then
these elements coincide graphically by induction hypothesis. Using the
assumptions~\eqref{eq:assumption}, the values of parameters fulfill the
following equalities
\begin{equation*}
    k = l, \Remainderm{p_k} = \Remainderm{q_l}, \Dividerm{p_1} + \Dividerm{p_0}
    + r_1 \n[p_0] = \Dividerm{q_1} + \Dividerm{q_0} + r_1 \n[p_0],
\end{equation*}
if $k = 1$, and
\begin{equation*}
\begin{gathered}
    k = l, \Remainderm{p_k} = \Remainderm{q_l}, \Dividerm{p_k} =
    \Dividerm{q_l}, p_{k - 1} = q_{l - 1}, \ldots, p_2 = q_2,\\
    p_{1} - 1 + \Dividerm{p_{0}} + r_1 \n[p_0] = q_{1} - 1 + \Dividerm{q_{0}} +
    r_1 \n[p_0],
\end{gathered}
\end{equation*}
otherwise. Adding the assumption $\Remainderm{p_0} = \Remainderm{q_0}$, the
sets of requirements are written in the following way
\begin{equation} \label{eq:conditions_k=1}
    k = l, p_1 + p_0 = q_1 + q_0,
\end{equation}
if $k = 1$, and
\begin{equation} \label{eq:conditions_k>1}
    k = l, p_k = q_l, p_{k - 1} = q_{l - 1}, \ldots, p_2 =
    q_2, m p_{1} + p_{0} = m q_{1} + q_{0},
\end{equation}
otherwise. If $p_0 = q_0$, then it follows from~\eqref{eq:conditions_k=1}
and~\eqref{eq:conditions_k>1} that the values of $p_i$ and $q_i$ coincide for
all $i = 0, 1, \ldots, k$, and elements $\s_1$ and $\s_2$ have the same normal
form.

\par Now let us assume that $p_0 \neq q_0$. As $\Remainderm{p_0} =
\Remainderm{q_0}$, then $r_2 \n[p_0] = r_2 \n[q_0]$ and it follows
from~\eqref{eq:assumption} that for any $u \in \InfiniteSet$ the equality holds
\begin{equation*}
    x_{\Remainderm{p_k}} \cdot \Restriction{\s_1}{x_{r_2 \n[p_0]}} \n[u] =
    x_{\Remainderm{q_k}} \cdot \Restriction{\s_2}{x_{r_2 \n[p_0]}} \n[u],
\end{equation*}
where the elements $\Restriction{\s_1}{x_{r_2 \n[p_0]}}$ and
$\Restriction{\s_2}{x_{r_2 \n[p_0]}}$ are defined by the
equalities~\eqref{eq:restriction_k=1_r_2} or~\eqref{eq:restriction_k>1_r_2}
depending on $k = 1$ or $k > 1$.

\par Let $k = 1$. By the equality at the line above the
elements
\begin{align*}
    \Restriction{\s_1}{x_{r_2 \n[p_0]}} & = f_1 ^{\Dividerm{p_1}} f_0 f_1
    ^{\Dividerm{p_0}} & &\text{and} & \Restriction{\s_2}{x_{r_2 \n[p_0]}} & =
    f_1 ^{\Dividerm{q_1}} f_0 f_1 ^{\Dividerm{q_0}}\\
\intertext{coincide graphically if and only if the elements}
    \s_1 & = f_1 ^{p_1} f_0 f_1 ^{p_0} & &\text{and} & \s_2 & = f_1 ^{q_1} f_0
    f_1 ^{q_0}
\end{align*}
coincide graphically.

\par Let $k > 1$, and let all integers $q_1, p_1, p_2, \ldots, p_{k - 1}$ are
divisible by $m$. Then $\Restriction{\s_1}{x_{r_2 \n[p_0]}}$ and
$\Restriction{\s_2}{x_{r_2 \n[p_0]}}$ can be written in the
form~\eqref{eq:normal_form}:
\begin{align*}
    \Restriction{\s_1}{x_{r_2 \n[p_0]}} & = f_1 ^{\Dividerm{p_k}} f_0 f_1 ^{m
    ^{k - 1} \Dividerm{p_{k - 1}} - 1} \ldots f_0 f_1 ^{m^2 \Dividerm{p_{2}} -
    1} f_0 f_1 ^{m \Dividerm{p_{1}} - 1} f_0 f_1 ^{\Dividerm{p_0}},\\
    \Restriction{\s_2}{x_{r_2 \n[p_0]}} & = f_1 ^{\Dividerm{p_k}} f_0 f_1 ^{m
    ^{k - 1} \Dividerm{p_{k - 1}} - 1} \ldots f_0 f_1 ^{m^2 \Dividerm{p_{2}} -
    1} f_0 f_1 ^{m \Dividerm{q_{1}} - 1} f_0 f_1 ^{\Dividerm{q_0}}.
\end{align*}
Similarly in this case the elements $\Restriction{\s_1}{x_{r_2 \n[p_0]}}$ and
$\Restriction{\s_2}{x_{r_2 \n[p_0]}}$ coincide graphically if and only if the
elements $\s_1$ and $\s_2$ coincide graphically.

\par Let $t$ be the maximal positive integer such that $\Remainderm{
\left[\frac{p _0}{m ^i} \right] } = \Remainderm{ \left[ \frac{q _0}{m ^i}
\right] }$ for all $i = 0, 1, \ldots, t - 1$, and in the case $k > 1$ all
integers $q_1, p_1, p_2, \ldots, p_{k - 1}$ are divisible by $m ^t$. As $p_0
\neq q_0$ then $t$ is a positive integer. Using the speculations above, it
follows from the assumption~\eqref{eq:assumption} that the following elements
define the same automaton transformations:
\begin{align*}
    \s_3 &= f_1 ^{\left[ \frac{p_1}{m ^t} \right]} f_0 f_1 ^{\left[
    \frac{p_0}{m ^t} \right]} & &\text{and} & \s_4 & = f_1 ^{\left[
    \frac{q_1}{m ^t} \right]} f_0 f_1 ^{\left[ \frac{q_0}{m ^t} \right]}
\end{align*}
if $k = 1$, and
\begin{align*}
    \s_5 & = f_1 ^{\left[ \frac{p_k}{m ^t} \right]} f_0 f_1 ^{m ^{k - 1} \left[
    \frac{p_{k - 1}}{m ^t} \right] - 1} \ldots f_0 f_1 ^{m^2 \left[
    \frac{p_2}{m ^t} \right] - 1} f_0 f_1 ^{m \left[ \frac{p_1}{m ^t} \right] -
    1} f_0 f_1 ^{\left[ \frac{p_0}{m ^t} \right]},\\
    \s_6 & = f_1 ^{\left[ \frac{p_k}{m ^t} \right]} f_0 f_1 ^{m ^{k - 1} \left[
    \frac{p_{k - 1}}{m ^t} \right] - 1} \ldots f_0 f_1 ^{m^2 \left[
    \frac{p_2}{m ^t} \right] - 1} f_0 f_1 ^{m \left[ \frac{q_1}{m ^t} \right] -
    1} f_0 f_1 ^{\left[ \frac{q_0}{m ^t} \right]},
\end{align*}
if $k > 1$. In addition, the elements $\s_1$ and $\s_2$ coincide graphically if
and only if $\s_3$ and $\s_4$ ($\s_5$ and $\s_6$, respectively) coincide
graphically .

\par As $t$ is maximal, then there are two possible cases:
\begin{enumerate}
\item $\Remainderm{ \left[\frac{p _0}{m ^t} \right] } \neq \Remainderm{ \left[
\frac{q _0}{m ^t} \right] }$,

\item $k > 1$, $\Remainderm{ \left[\frac{p _0}{m ^t} \right] } = \Remainderm{
\left[ \frac{q _0}{m ^t} \right] }$, and one of $q_1, p_1, p_2, \ldots, p_{k -
1}$ is not divisible by $m ^{t + 1}$.
\end{enumerate}

It follows from Item 1 that in the case $\Remainderm{ \left[\frac{p _0}{m ^t}
\right] } \neq \Remainderm{ \left[ \frac{q _0}{m ^t} \right] }$ the
contradiction with the assumption~\eqref{eq:assumption} follows from the
equality $\s_3 = \s_4$ (or $\s_5 = \s_6$).

\par In the second case $k > 1$ and $\Remainderm{ \left[\frac{p _0}{m ^t}
\right] } = \Remainderm{ \left[ \frac{q _0}{m ^t} \right] }$. Let us consider
the input word $v = x_{r_2 \n[ {\left[\frac{p _0}{m ^t} \right]} ]}$. At least
one of the elements $\Restriction{\s_5}{v}$ and $\Restriction{\s_6}{v}$ is
reducible, and have normal form of length $\left( k - 1 \right)$. If another
element is irreducible, the contradiction with~\eqref{eq:assumption} follows
from the induction hypothesis. Hence both elements $\Restriction{\s_5}{v}$ and
$\Restriction{\s_6}{v}$ are reducible. It follows from the note
on~\eqref{eq:restriction_k>1_r_2} that normal forms of these elements end with
$f_0 f_1 ^{\left[ \frac{p_0}{m ^{t + 1}} \right]}$ and $f_0 f_1 ^{\left[
\frac{q_0}{m ^{t + 1}} \right]}$, respectively. Then it follows from the
induction hypothesis that the equality
\[
    \left[ \frac{p_0}{m ^{t + 1}} \right] = \left[ \frac{q_0}{m ^{t + 1}}
    \right]
\]
holds. Moreover, by assumptions, the equalities $\Remainderm{ \left[\frac{p
_0}{m ^i} \right] } = \Remainderm{ \left[ \frac{q _0}{m ^i} \right] }$ hold for
all $i = 0, 1, \ldots, t$, whence $p_0 = q_0$. Combined with the
requirements~\eqref{eq:conditions_k>1}, we have the set of requirements
\[
    k = l, p_k = q_l, p_{k - 1} = q_{l - 1}, \ldots, p_2 = q_2, p_{1} = q_{1},
    p_{0} = q_{0},
\]
i.e. the normal forms of $\s_1$ and $\s_2$ coincide graphically.

\par The proposition is completely proved.
\end{proof}

\subsection{Cayley graph} \label{subsect:semigroup_Cayley_graph}

\par In this subsection we construct the Cayley graph
$G_{\Semigroup{\Intermediate}}$ of the semigroup $\Semigroup{\Intermediate}$.
$\Semigroup{\Intermediate}$ is a monoid, and the root of graph is the identity,
that belongs to the semigroup. As we apply automaton transformation from right
to left, then we will read the labels of path in the same order. For example,
the edges labelled by $f_1 \text{\ndash} f_1 \text{\ndash} f_0$ denote the path
$f_0 f_1 ^2$. It follows from Proposition~\ref{prop:normal_form} and
Proposition~\ref{prop:normality_of_form} that an arbitrary element $\s \in
\Semigroup{\Intermediate}$ can be unambiguously reduced to the
form~\eqref{eq:normal_form}. Hence any path
without loops should define the semigroup element in normal form.

\par The graph $G_{\Semigroup{\Intermediate}}$ consists of subgraphs $E_i$, $i \ge
0$. An arbitrary path in $G_{\Semigroup{\Intermediate}}$ walks through groups
of $E_i$, $i = 0, 1, \ldots$, connected by edges labelled $f_0$, and each group
consists of several copies of $E_i$, connected by edges labelled $f_1$. The
path, defined by $p_i$ copies of $E_i$, corresponds to the subword $f_1 ^{m ^i
p_i - 1}$ in the semigroup word written in the form~\eqref{eq:normal_form}.

\par
\begin{figure}[t]
  \centering
  \includegraphics*{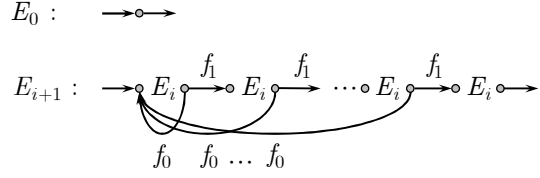}
  \caption{The graphs $E_i$, $i \ge 0$}
  \label{fig:graph_blocks}
\end{figure}

 The structure of the graphs $E_i$, $i \ge 0$, is shown on
Figure~\ref{fig:graph_blocks}. The rightmost and the leftmost arrows on the
figure don't belong to $E_i$ and denote edges, that enter and output from the
graph $E_i$. The shaded circles before and after the graph $E_i$ denote the
rightmost and the leftmost vertices of $E_i$. The graph $E_0$ includes a unique
vertex, and does not have edges. The graph $E_{i + 1}$ is constructed as $m$
copies of $E_i$, that are sequentially connected by edges labelled by $f_1$,
and the rightmost vertex of each of the first $\left( m - 1 \right)$ graphs
$E_i$ is connected to the leftmost vertex of the first graph $E_i$ by the edge
labelled by $f_0$.

\begin{lemma} \label{lem:Ei_f1 count}
For all $i \ge 0$ the graph $E_i$ includes $m ^i - 1$ edges labelled by the
symbol $f_1$.
\end{lemma}

\begin{proof}
We prove the lemma by induction on $i$. Clearly $E_0$ includes $0 = m^0 - 1$
edges labelled by the symbol $f_1$. For $i \ge 0$ the graph $E_i$ includes
\[
    \underbrace{m \left( m^{i - 1} - 1 \right)} _{\text{$m$ copies of $E_{i -
    1}$}} + \underbrace{m - 1} _{\text{edges between $E_{i - 1}$s}} = m ^i - 1
\]
edges labelled by $f_1$. The lemma is proved.
\end{proof}

\par  Thus each graph $E_i$ can be presented as a direct path with $\left( m
^i - 1 \right)$ edges labelled by $f_1$; and the edge labelled by $f_0$ outputs
from each vertex, excepting the rightmost, and enters one of the previous
vertices.

\begin{lemma} \label{lem:Ei_reducing}
    Let $i$, $p$ be arbitrary integers such that $i > 0$ and $0 \le p < m^i -
    1$. Then the path $P = f_0 f_1 ^p$ in the graph $E_i$ that starts from the
    leftmost vertex, includes the loop $L = f_0 f_1 ^{t_2 \n[p + 1] - 1}$.
\end{lemma}

\begin{proof}
As $p < m^i - 1$ then the path $f_0 f_1 ^p$ belongs to $E_i$ and doesn't
include the rightmost vertex of $E_i$. It follows from the note after
Lemma~\ref{lem:Ei_f1 count} that an arbitrary edge labelled by $f_0$ forms a
loop in the graph $E_i$.

\par We prove the second statement by induction on $i$. For $i = 1$ the path is
$f_0 f_1 ^p$, where $0 \le p < m - 1$. Using definitions of $t_1$ and $t_2$,
the following equalities hold
\begin{align*}
    t_1 \n[p + 1] & = 0, & t_2 \n[p + 1] & = p + 1,
\end{align*}
whence $t_2 \n[p + 1] - 1 = p$. It follows from Figure~\ref{fig:graph_blocks}
that $P$ is a loop, and, consequently, $L = P = f_0 f_1 ^p$.

\par Let $i > 1$. The graph $E_i$ consists of $m$ copies of $E_{i - 1}$, and
there are two possible cases for the edge labelled by $f_0$: it is contained
within one of $E_{i - 1}$ or connects the rightmost vertex of one of $E_{i -
1}$ with the leftmost vertex of $E_i$.

\par Each $E_{i - 1}$ includes $\n[m ^{i - 1} - 1]$ edges labelled by $f_1$,
and in the first case the equality
\[
    p = q \cdot m ^{i - 1} + r
\]
holds for $0 \le q < m$, $0 \le r < m ^{i - 1} - 1$. Let $P' = f_0 f_1 ^{r}$ be
the path that starts from the leftmost vertex of $\n[q + 1]$-th copy of $E_{i -
1}$. By induction hypothesis $L = f_0 f_1 ^{t_2 \n[r + 1] - 1}$. As $\n[r + 1]
< m ^{i - 1}$ and $\n[r + 1]$ is not divisible by $m ^{i - 1}$, then
\[
    t_1 \n[p + 1] = t_1 \n[q \cdot m ^{i - 1} + {\n[r + 1]}] = t_1 \n[r + 1]
    \le i - 2.
\]
Therefore the equalities hold
\begin{multline*}
    t_2 \n[p + 1] = \n[q \cdot m ^{i - 1} + {\n[r + 1]}] \mod m ^{t_1 \n[r + 1]
    + 1} \\
    = \n[r + 1] \mod m ^{t_1 \n[r + 1] + 1} = t_2 \n[r + 1].
\end{multline*}
Hence, the path $P$ includes the loop $L = f_0 f_1 ^{t_2 \n[r + 1] - 1} = f_0
f_1 ^{t_2 \n[p + 1] - 1}$.

\par In the second case the path $P$ is a loop. Similarly, the equality
\[
    p = q \cdot m ^{i - 1} + \n[m ^{i - 1} - 1]
\]
holds for $0 \le q < m - 1$. As $p + 1 = \n[q + 1] m ^{i - 1}$, we have
\begin{align*}
    t_1 \n[p + 1] & = i - 1 & &\text{and} & t_2 \n[p + 1] & = \n[q + 1] m ^{i
    - 1} = p + 1.
\end{align*}
Therefore $L = f_0 f_1 ^{t_2 \n[p + 1] - 1} = f_0 f_1 ^p = P$.
\end{proof}

\par
\begin{figure}[t]
  \centering
  \includegraphics*{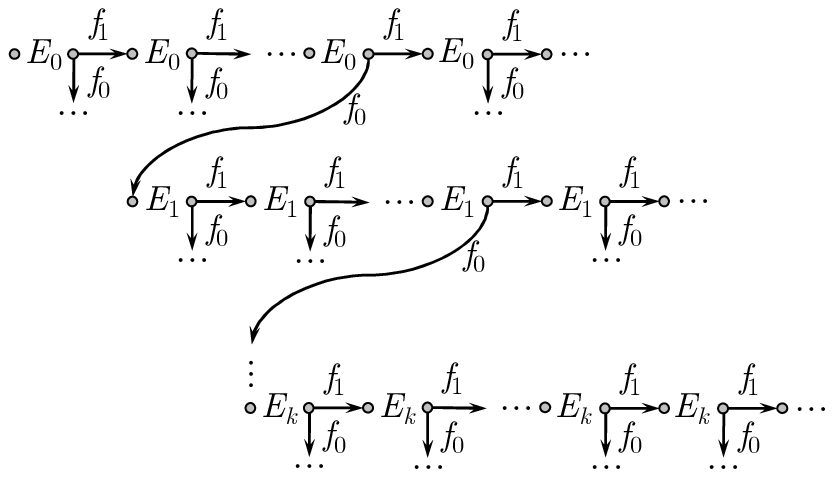}
  \caption{The graph $G_{\Semigroup{\Intermediate}}$}
  \label{fig:graph_Cayley}
\end{figure}

 The Cayley graph $G_{\Semigroup{\Intermediate}}$ is shown on
Figure~\ref{fig:graph_Cayley}. The generator $e$ gives loops labelled by $e$ on
each vertex, and we don't show these edges. The graph
$G_{\Semigroup{\Intermediate}}$ can be conditionally separated into lines,
where $i$-th line, $i \ge 0$, consists of copies of $E_i$. These graphs are
connected by edges labelled by $f_1$, and the edges labelled by $f_0$ allow to
pass to the next line. The leftmost vertex of zero line is the root of
$G_{\Semigroup{\Intermediate}}$ and corresponds to the semigroup identity.

\begin{proposition} \label{prop:Cayley_graph}
    Let $P$ be an arbitrary path in $G_{\Semigroup{\Intermediate}}$ such that
    it starts from the root vertex, and let it denotes the semigroup word $\s$.
    Then $P$ includes the path $P'$ without loops such that it denotes the
    semigroup element $\s'$ written in the normal form~\eqref{eq:normal_form}
    and is equivalent to $\s$.
\end{proposition}

\begin{proof}
At first, we show that the path $P$ without loops denotes the semigroup element
written in the normal form~\eqref{eq:normal_form}. It follows from
Lemma~\ref{lem:Ei_reducing} and structure of $E_i$ that loops are created only
by edges labelled by $f_0$ that are located within $E_i$. If $P$ doesn't
include the edges labelled by $f_0$, then it's located at zero line of
$G_{\Semigroup{\Intermediate}}$, and $P$ denotes the semigroup word $\s = f_1
^p$ for some $p \ge 0$. Clearly $\s$ is written in the normal form.

\par Otherwise, let $P$ ends at the $k$-th line of
$G_{\Semigroup{\Intermediate}}$, $k \ge 1$. Then all edges labelled by $f_0$,
that belong to $P$, connect the lines of $G_{\Semigroup{\Intermediate}}$. Let
the path $P$ goes through $p_i$ copies of $E_i$ at $i$-th line, $0 \le i < k$.
It follows from Lemma~\ref{lem:Ei_f1 count} that the subpath of $P$ at $i$-th
line denotes the semigroup word $f_1 ^{p_i m ^i - 1}$. Therefore, the path
denotes the semigroup word $\s = f_1 ^{p_k} f_0 f_1 ^{m ^{k - 1} p_{k - 1} - 1}
f_0 \ldots f_1 ^{m ^2 p_2 - 1} f_0 f_1 ^{m p_1 - 1} f_0 f_1 ^{p_0}$, where $p_0
\ge 0$, $p_i \ge 1$, $1 \le i < k$, and $p_k \ge 0$ is the count of edges of
$P$ at $k$-th line. Similarly, $\s$ is written in the
form~\eqref{eq:normal_form}.

\par Let $\s$ be an arbitrary word over the alphabet $\left\{ f_0, f_1
\right\}$:
\begin{equation*}
    \s = f_1 ^{p_{2k}} f_0 ^{p_{2k - 1}} f_1 ^{p_{2k - 2}} \ldots f_0
    ^{p_{1}} f_1 ^{p_{0}},
\end{equation*}
where $k \ge 0$, $p_0, p_{2k} \ge 0$, $p_i \ge 1$, $i = 1, 2, \ldots, {2k -
1}$, and let us consider the path $P$ that denotes $\s$. We show that
sequential cancelling of loops in $P$ coincide with the executing of
Algorithm~\ref{alg:reducing}. The variable $j$ denotes the current line of
$G_{\Semigroup{\Intermediate}}$, and $r$ denotes the length of subpath at
$j$-th line. As the leftmost vertex of $i$-th line, $i > 0$, is located within
$E_i$ and has the loop labelled by $f_0$, then this loop is removed by
operations at the
lines~\ref{alg_reducing:if_odd_begin}--\ref{alg_reducing:if_odd_end} of the
algorithm. The check at the line~\ref{alg_reducing:if_even1} means that subpath
at $j$-th line reaches the rightmost vertex of $E_i$. It follows from
Lemma~\ref{lem:Ei_reducing} that actions of the
lines~\ref{alg_reducing:if_even2_begin}--\ref{alg_reducing:if_even2_end} are
realized by reduction of loops $f_0 f_1 ^{p_i}$ inside of the graph $E_i$.
Hence, $P$ includes the path $P'$ without loops, that denotes the semigroup
word $\s'$ in the form~\eqref{eq:normal_form} that is equivalent to $\s$.
\end{proof}

\subsection{Proof of Theorem~\ref{th:semigroup}} \label{subsect:semigroup_proof}

\begin{proof}[Proof of Theorem~\ref{th:semigroup}]
From Lemma~\ref{prop:all_relations} follows that in the semigroup
$\Semigroup{\Intermediate}$ the relations
\[
    r \left( k, p_{k + 2}, p_{k + 1}, p_{k}, \ldots, p_{1} \right)
\]
where $k \ge 0$, $1 \le p_{k + 2} \le {m - 1}$, $p_{k + 1} \ge 0$, $p_i \ge 1$,
$i = 1, 2, \ldots, {k - 1}$, hold. In Proposition~\ref{prop:normal_form} it is
shown that, using these relations, each element can be reduced to the
form~\eqref{eq:normal_form}. On the other hand, it is proved in
Proposition~\ref{prop:normality_of_form} that two semigroup elements written in
the form \eqref{eq:normal_form} define the same automaton transformation over
the set $\InfiniteSet$ if and only if they coincide graphically. Hence, the
form~\eqref{eq:normal_form} is the normal form of elements of
$\Semigroup{\Intermediate}$, and each semigroup element can be unambiguously
reduced to the form~\eqref{eq:normal_form}.

\par The set of relations~\eqref{eq:all_relations} is not minimal. It is proved
in Proposition~\ref{prop:defining_relations} and Lemma~\ref{prop:all_relations}
that in the semigroup $\Semigroup{\Intermediate}$ the
relations~\eqref{eq:all_relations} may be derived from the
set~\eqref{eq:defining_relations} of relations:
\begin{align*}
    R_A & \Pair{k}{p}, k \ge 0, p = 1, 2, \ldots, {m - 1}, & R_B &
    \Argument{k}, k \ge 0.
\end{align*}
The structure of the Cayley graph of the semigroup~$\Semigroup{\Intermediate}$
is considered in Subsection~\ref{subsect:semigroup_Cayley_graph}. It follows
from Figure~\ref{fig:graph_blocks} and Figure~\ref{fig:graph_Cayley} that the
edges, that realize the reducing of semigroup words, belong to the graphs
$E_k$, $k \ge 0$. Each relation of type $A$ of length $k$ substantiates the
edge labelled by $f_0$ that forms a loop in the graph $E_{k + 1}$. The
relations of type $B$ of length $k$ allow to connect the graphs $E_{k + 1}$ at
the $\left( k + 1 \right)$-th line of the graph
$G_{\Semigroup{\Intermediate}}$, one relation per line. Therefore, the
set~\eqref{eq:defining_relations} is minimal, that is no one relation follows
from the others. Thus, the infinite set of
relations~\eqref{eq:defining_relations} is the system of defining relations,
and the semigroup $\Semigroup{\Intermediate}$ is infinitely presented.

\par The automaton transformation $e$ is the identity, whence
$\Semigroup{\Intermediate}$ is an infinitely presented
monoid.

\par To solve the word problem in $\Semigroup{\Intermediate}$, it is
necessary to reduce semigroup words $\s_1$ and $\s_2$ to normal
form~\eqref{eq:normal_form}, and then to check them for graphical equality.
From Proposition~\ref{prop:normal_form} follows that count of steps, required
by both reductions, is equivalent
\[
    \mathcal{O} \Argument{\Argument{\left| \s_1 \right| + \left| \s_1 \right|}
    \log_m \Argument{\left| \s_1 \right| + \left| \s_1 \right|} },
\]
and the word problem is solved in no more than quadratic time.
\end{proof} 

\begin{proof}[Proof of Corollary~\ref{cor:rewriting_system}]
It follows from Proposition~\ref{prop:normal_form} and
Algorithm~\ref{alg:reducing} that an arbitrary semigroup element $\s$ can be
reduced to the normal form~\eqref{eq:normal_form} by applying the
relations~\eqref{eq:all_relations}. On the other hand, the element $\s$ is
written in the normal form if and only if it doesn't include left part of any
relation~\eqref{eq:all_relations}. Therefore, the set of relations
\[
    r \left( k, p_{k + 2}, p_{k + 1}, p_{k}, \ldots, p_{1} \right)
\]
for all possible values $k \ge 0$, $1 \le p_{k + 2} \le {m - 1}$, $p_{k + 1}
\ge 0$, $p_i \ge 1$, $i = 1, 2, \ldots, {k - 1}$, is the rewriting system of
the monoid $\Semigroup{\Intermediate}$. It follows from
Theorem~\ref{th:semigroup} that elements of the form~\eqref{eq:normal_form} is
in bijection with elements of $\Semigroup{\Intermediate}$, whence this
rewriting system is complete.
\end{proof}

\section{Growth of $\Intermediate$ and $\Semigroup{\Intermediate}$}
\label{sect:growth_functions}

We derive, in this section, the growth series of the semigroup
$\Semigroup{\Intermediate}$ and the automaton $\Intermediate$, as well as the
asymptotics of the growth functions $\GrowthSemigroup{\Intermediate}$ and
$\Growth{\Intermediate}$.

\subsection{Growth series} \label{subsect:growth_series}

\par Natural system of generators of the monoid $\Semigroup{\Intermediate}$
includes the identity, and it follows from
Proposition~\ref{prop:monoid_growth_order} that the equality holds
\begin{equation*}
    \Growth{\Intermediate} \n = \GrowthSemigroup{\Intermediate} \n, \quad
    \PositiveInteger.
\end{equation*}
Obviously, it implies that $\Gamma _{\Intermediate} \Argument{X} = \Gamma
_{\Semigroup{\Intermediate}} \Argument{X}$.

\par At first, we derive the growth series $\Delta _{\Semigroup{\Intermediate}}
\Argument{X}$ for the word growth function of~$\Semigroup{\Intermediate}$. It
follows from Theorem~\ref{th:semigroup} that each semigroup element $\s$ can be
unambiguously reduced to the form~\eqref{eq:normal_form}. We arrange all
semigroup elements by length of their normal form, and the growth series that
count elements of length $l$, $l \ge 0$, are listed in the following table:
\begin{align*}
    l &= 0: & & f_1 ^{p_0} & \frac{1}{1 - X} & ;\\
    l &= 1: & & f_1 ^{p_1} \underbrace{f_0 f_1 ^{p_0}} _{p_0 + 1} & \frac{1}{1 -
    X} & \cdot \frac{X}{1 - X};\\
    l &= 2: & & f_1 ^{p_2} \underbrace{f_0 f_1 ^{m p_1 - 1}} _{m p_1}
    \underbrace{f_0 f_1 ^{p_0}} _{p_0 + 1} & \frac{1}{1 - X} & \cdot \frac{X}{1
    - X} \cdot \frac{X ^m}{1 - X ^m};\\ & \ldots & & \ldots & & \ldots\\
    l &= k: & & {\begin{array}{*{20}l}
        {f_1 ^{p_{k}} \underbrace{f_0 f_1 ^{m ^{k - 1} p_{k - 1} - 1}} _{m ^{k
        - 1} p_{k - 1}} \ldots }\\
        {\qquad \quad \underbrace{f_0 f_1 ^{m ^2 p_{2} - 1}} _{m ^2 p_{2}}
        \underbrace{f_0 f_1 ^{m p_{1} - 1}} _{m p_{1}} \underbrace{f_0 f_1
        ^{p_0}} _{p_0 + 1}}
    \end{array}}&
    \frac{1}{1 - X} & \cdot \prod \limits_{i = 0}^{k - 1} \frac{X ^{m ^i}}{1 -
    X ^{m ^i}} ;\\
    & \ldots & & \ldots & & \ldots\\
\end{align*}
First column includes the length $l$ of normal form, second --- the set of the
semigroup elements in the form~\eqref{eq:normal_form} of length $l$, and the
corresponding growth series are listed in third column. Let $\s$ be an
arbitrary semigroup word in the form~\eqref{eq:normal_form}:
\[
    f_1 ^{p_k} f_0 f_1 ^{m ^{k - 1} p_{k - 1} - 1} \ldots f_0 f_1 ^{m ^{i}
    p_{i} - 1} \ldots f_0 f_1 ^{m ^2 p_{2} - 1} f_0 f_1 ^{m p_{1} - 1} f_0 f_1
    ^{p_0}
\]
where $k \ge 0$, $p_0 \ge 0$, $p_k \ge 0$, $p_i \ge 1$, $i = 1, 2, \ldots, {k -
1}$. Every subword $f_0 f_1 ^{m ^i p_{i} - 1}$, $1 \le i \le {k - 1}$, $p_{i}
\ge 1$, has length $m ^{i} p_{i}$ and is counted by the growth series
\[
    \frac{X ^{m ^i}}{1 - X ^{m ^i}};
\]
The end $f_0 f_1 ^{p_0}$, $p_0 \ge 0$, has the length $\Argument{p_0 + 1}$ and
is counted by the growth series $\frac{X}{1 - X}$, and the beginning $f_1
^{p_k}$, $p_k \ge 0$, is counted by $\frac{1}{1 - X}$. Hence, the word growth
series of $\Semigroup{\Intermediate}$ is
\begin{multline*}
    \Delta _{\Semigroup{\Intermediate}} \Argument{X} = \sum \limits_{k \ge
    0} {\frac{1}{1 - X} \prod \limits_{i = 0}^{k - 1} \frac{X^{m^i}}{1 -
    X^{m^i}}} = \frac{1}{1 - X} \left( 1 + \frac{X}{1 - X} \left( 1 +
    \frac{X^m}{1 - X ^m} \cdot \right. \right.\\
    \left. \left. \left( 1 + \frac{X ^{m^2}}{1 - X ^{m^2}} \left( 1 +
    \frac{X ^{m^3}}{1 - X ^{m^3}} \left(1 + \ldots \right) \right) \right)
    \right) \right),
\end{multline*}
that proves Corollary~\ref{cor:word_generating_function}.

\par It follows from the note at the beginning of this subsection, that
\begin{multline*}
    \Gamma _{\Intermediate} \Argument{X} = \Gamma _{\Semigroup{\Intermediate}}
    \Argument{X} = \frac{1}{1 - X} \Delta _{\Semigroup{\Intermediate}}
    \Argument{X} = \frac{1}{\Argument{1 - X} ^2} \left( 1 + \frac{X}{1 - X}
    \cdot \right.\\
    \left. \left( 1 + \frac{X ^m}{1 - X ^m} \left( 1 + \frac{X ^{m^2}}{1
    - X ^{m^2}} \left( 1 + \frac{X ^{m^3}}{1 - X ^{m^3}} \left(1 + \ldots
    \right) \right) \right) \right) \right),
\end{multline*}
that completes the proof of Theorem~\ref{th:generating_functions}.

\par Let us denote the series
\[
    1 + \frac{X}{1 - X} \left( 1 + \frac{X ^m}{1 - X ^m} \left( 1 + \frac{X ^{m
    ^2}}{1 - X ^{m ^2}} \left( 1 + \frac{X ^{m ^3}}{1 - X ^{m ^3}} \left( 1 +
    \ldots \right) \right) \right) \right)
\]
by the symbol $S \n[X]$, and the growth series are defined by the following
equalities
\begin{align*}
    \Delta _{\Semigroup{\Intermediate}} \Argument{X} & = \frac{1}{1 - X} S
    \n[X], &
    \Gamma _{\Intermediate} \Argument{X} & = \Gamma
    _{\Semigroup{\Intermediate}} \Argument{X} = \frac{1}{\left( 1 - X \right)
    ^2} S \n[X].
\end{align*}

\subsection{Proof of Corollary~\ref{cor:numerical_properties}}
\label{subsect:growth_proof_corollary}

\par Let $g : \Natural \to \Natural$ be an arbitrary function of a natural
argument, and $G \n[X] = \sum \limits _{n \ge 0} g \n X^n$ be its generating
function. Then the following equalities hold
\begin{subequations}
\begin{equation} \label{eq:generating_functions_n+1}
    \sum \limits _{n \ge 0} g \n[n + 1] X^n = \sum \limits _{n \ge 1} g \n
    X^{n - 1} = \frac{1}{X} \left( G \n[X] - g \n[0] \right),
\end{equation}
and
\begin{multline} \label{eq:generating_functions_n/m}
    \sum \limits _{n \ge 0} g \Argument{\Dividerm{n}} X^n = \sum \limits
    _{n \ge 0} g \n \left( X ^{mn} + X ^{mn + 1} + \ldots + X ^{mn + m - 1}
    \right)\\
    = \left( 1 + X + \ldots + X ^{m - 1} \right) \sum \limits _{n \ge 0} g \n
    X ^{mn} = \frac{1 - X ^m}{1 - X} \, G \n[X ^m].
\end{multline}
Let $p$ be a positive integer, $r = \exp \Argument{\frac{2 \pi i}{p}}$ be a
primary $p$-th root of the identity. Applying the method of power series
multisection~\cite{Riordan1968}, for any $0 \le k < p$ the $k$-th section of
the series $G \n[X]$ is defined by the equality
\[
    \sum \limits _{n \ge 0} g \n[k + np] X ^{k + np} = \frac{1}{p} \sum \limits
    _{j = 1} ^{p} r ^{p - kj} G \n[r ^j X].
\]
Hence the generating series of $g' \n = g \n[mn]$, $\PositiveInteger$, is
described by the following equality
\begin{equation} \label{eq:generating_functions_mn}
    \sum \limits _{n \ge 0} g \n[m n] X ^n = \frac{1}{m} \sum \limits _{j = 1}
    ^{m} r ^{m} G \n[r ^j X ^{\frac{1}{m}}] = \frac{1}{m} \sum \limits _{j = 1}
    ^{m} G \n[r ^j X ^{\frac{1}{m}}],
\end{equation}
\end{subequations}
where $r = \exp \Argument{\frac{2 \pi i}{m}}$.

\begin{proof}[Proof of Corollary~\ref{cor:numerical_properties}]
\par \ref{cor_item:functional_equation}. Let us write the
equality~\eqref{eq:first_difference_aymptotics} in terms of generating
functions, and it's enough to prove, that the equality
\begin{equation} \label{eq:generating_functions_2}
    {\sum \limits _{n \ge 0} \WordGrowthSemigroup{\Intermediate}
    \n[\Dividerm{n}] X^n} + {\sum \limits _{n \ge 0}
    \WordGrowthSemigroup{\Intermediate} \n X^n} - {\sum \limits _{n \ge 0}
    \WordGrowthSemigroup{\Intermediate} \n[n + 1] X^n} = 0
\end{equation}
holds. It is proved in Theorem~\ref{th:generating_functions} that $\Delta
_{\Semigroup{\Intermediate}} \Argument{X} = \sum \limits_{n \ge 0}
\WordGrowthSemigroup{\Intermediate} \n X ^n = \frac{1}{1 - X} S \n[X]$, whence
the following equality holds
\begin{multline*}
    \Delta _{\Semigroup{\Intermediate}} \Argument{X ^m} = \frac{1}{1 - X ^m}
    S \n[X ^m] \\
    = \frac{1}{1 - X ^m} \left(1 + \frac{X ^m}{1 - X ^m} \left( 1 + \frac{X ^{m
    ^2}}{1 - X ^{m ^2}} \left( 1 + \frac{X ^{m ^3}}{1 - X ^{m ^3}} \left( 1 +
    \ldots \right) \right) \right) \right) \\
    = \frac{1}{1 - X ^m} \cdot \frac{1 - X}{X} \cdot \left( S \n[X] - 1
    \right).
\end{multline*}
Using the equalities~\eqref{eq:generating_functions_n+1},
\eqref{eq:generating_functions_n/m} and the equality at the line above, we have
\begin{multline*}
    \frac{1 - X ^m}{1 - X} \Delta _{\Semigroup{\Intermediate}} \n[X ^m] +
    \Delta _{\Semigroup{\Intermediate}} \n[X] - \frac{1}{X} \left( \Delta
    _{\Semigroup{\Intermediate}} \n[X] - 1 \right) \\
    = \frac{1}{X} \left( S \n[X] - 1 \right) + \frac{1}{1 - X} S \n[X] \left( 1
    - \frac{1}{X} \right) + \frac{1}{X} = 0.
\end{multline*}
Hence, the equality~\eqref{eq:generating_functions_2} holds, and the statement
of the item is true.

\par \ref{cor_item:gamma_one_half}. Similarly to the
previous item, we prove that the following equality holds for the power series
\begin{equation} \label{eq:gamma_one_half_series}
    \sum \limits _{n \ge 0} \left( {m \Growth{\Intermediate} \n + 1 } \right) X
    ^n = \sum \limits _{n \ge 0} \WordGrowthSemigroup{\Intermediate} \n[m
    \Argument{n + 1}] X ^n.
\end{equation}
The left-side series of~\eqref{eq:gamma_one_half_series} can be easily
calculated:
\[
    \sum \limits _{n \ge 0} \left( {m \Growth{\Intermediate} \n + 1 } \right) X
    ^n = m \sum \limits _{n \ge 0} {\Growth{\Intermediate} \n X ^n} +
    \frac{1}{1 - X} = \frac{m}{\left( 1 - X \right) ^2} S \n[X] + \frac{1}{1 -
    X}.
\]

\par Let $r$ be a primary $m$-th root of $1$. Then the following equality
holds
\[
    \Argument{r ^j X ^{\frac{1}{m}}} ^{m ^k} = r ^{j m ^k} X ^{\frac{1}{m} m
    ^k} = X ^{m ^{k - 1}},
\]
for any $1 \le j \le m$, $k \ge 1$, whence we have
\begin{multline*}
    \Delta _{\Semigroup{\Intermediate}} \Argument{r ^j X ^{\frac{1}{m}}} =
    \frac{1}{1 - \Argument{r ^j X ^{\frac{1}{m}}}} S \n[\Argument{r ^j X
    ^{\frac{1}{m}}}] = \frac{1}{1 - \Argument{r ^j X ^{\frac{1}{m}}}}\\
    \cdot \left(1 + \frac{\Argument{r
    ^j X ^{\frac{1}{m}}}}{1 - \Argument{r ^j X ^{\frac{1}{m}}}} \left( 1 +
    \frac{X}{1 - X} \left( 1 + \frac{X ^{m ^2}}{1 - X ^{m ^2}}
    \left( 1 + \ldots \right) \right) \right) \right)\\
    = \frac{1}{1 - \Argument{r ^j X ^{\frac{1}{m}}}} + \frac{\Argument{r ^j X
    ^{\frac{1}{m}}}}{\left( 1 - \Argument{r ^j X ^{\frac{1}{m}}} \right) ^2} S
    \n[X].
\end{multline*}
\par Let us consider two power series
\begin{align*}
    A \n[X] & = \sum \limits _{n \ge 0} X ^n = \frac{1}{1 - X};&
    B \n[X] & = \sum \limits _{n \ge 0} n X ^n = \frac{X}{\left( 1 - X \right)
    ^2}.
\end{align*}
It follows from~\eqref{eq:generating_functions_mn} that the following
equalities hold
\begin{equation*}
    \frac{1}{m} \sum \limits _{j = 1} ^{m} \left( \frac{1}{1 - \Argument{r ^j X
    ^{\frac{1}{m}}}} \right) = \frac{1}{m} \sum \limits _{j = 1} ^{m} A \n[r
    ^j X ^{\frac{1}{m}}] = \sum \limits _{n \ge 0} a \n[nm] X ^n = \frac{1}{1 -
    X},
\end{equation*}
and
\begin{multline*}
    \frac{1}{m} \sum \limits _{j = 1} ^{m} {\frac{\Argument{r ^j X
    ^{\frac{1}{m}}}}{\left( 1 - \Argument{r ^j X ^{\frac{1}{m}}} \right) ^2}}
    = \frac{1}{m} \sum \limits _{j = 1} ^{m} B \n[r ^j X ^{\frac{1}{m}}]\\
    = \sum \limits _{n \ge 0} b \n[nm] X ^n = \sum \limits _{n \ge 0} nm X ^n =
    \frac{mX}{\left( 1 - X \right) ^2}.
\end{multline*}
Applying the equality~\eqref{eq:generating_functions_mn} and the equalities
proved above, we may write out the growth series of the function
$\WordGrowthSemigroup{\Intermediate} \n[m n]$:
\begin{multline*}
    \sum \limits _{n \ge 0} \WordGrowthSemigroup{\Intermediate} \n[m n] X ^n =
    \frac{1}{m} \sum \limits _{j = 1} ^{m} \Delta _{\Semigroup{\Intermediate}}
    \n[r ^j X ^{\frac{1}{m}}] = \frac{1}{m} \sum \limits _{j = 1} ^{m} \left(
    \frac{1}{1 - \Argument{r ^j X ^{\frac{1}{m}}}} \right) + \\
    + \frac{1}{m} \sum \limits _{j = 1} ^{m} \left( \frac{\Argument{r ^j X
    ^{\frac{1}{m}}}}{\left( 1 - \Argument{r ^j X ^{\frac{1}{m}}} \right) ^2}
    \right) S \n[X] = \frac{1}{1 - X} + \frac{mX}{\left( 1 - X \right) ^2} S
    \n[X].
\end{multline*}
Using the equality at the line above and~\eqref{eq:generating_functions_n+1},
we have
\begin{multline*}
    \sum \limits _{n \ge 0} \WordGrowthSemigroup{\Intermediate} \n[m \left( n +
    1 \right)] X ^n = \frac{1}{X} \left( \frac{1}{1 - X} + \frac{mX}{\left( 1 -
    X \right) ^2} S \n[X] - 1 \right)\\
    = \frac{1}{1 - X} + \frac{m}{\left( 1 - X \right) ^2} S \n[X].
\end{multline*}
Thus the left-hand and right-hand series of~\eqref{eq:gamma_one_half_series}
coincide, and the equality~\eqref{eq:gamma_one_half_series} is true, whence the
statement of Item~\ref{cor_item:gamma_one_half} is true.

\par \ref{cor_item:second_difference_partitions}. It follows
from~\eqref{eq:growths_semigroup}, that the second finite difference of
$\GrowthSemigroup{\Intermediate}$ is the first finite difference of
$\WordGrowthSemigroup{\Intermediate}$, i.e.
\[
    \GrowthSemigroup{\Intermediate} ^{\n[2]} \n =
    \WordGrowthSemigroup{\Intermediate} \n -
    \WordGrowthSemigroup{\Intermediate} \n[n - 1],
\]
for all $n \ge 1$, and let us assume $\GrowthSemigroup{\Intermediate} ^{\n[2]}
\n[0] = 1$. Let us denote the generating function of the function
$\Growth{\Intermediate} ^{\n[2]}$ by the symbol $\Gamma ^{\n[2]} \n[X]$. As
$\Growth{\Intermediate} ^{\n[2]} \n[0] = \WordGrowthSemigroup{\Intermediate}
\n[0] = 1$, then the following equality holds
\begin{equation*}
    \Gamma ^{\n[2]} \n[X] = \left( {1 - X} \right) \Delta
    _{\Semigroup{\Intermediate}} \n[X],
\end{equation*}
whence $\Gamma ^{\n[2]}$ can be presented as infinite sum of finite products:
\begin{multline} \label{eq:second_difference_partitions}
    \Gamma ^{\n[2]} \n[X] = 1 + \frac{X}{1 - X} \left( 1 + \frac{X ^m}{1 - X ^m}
    \left( 1 + \frac{X ^{m ^2}}{1 - X ^{m ^2}} \right. \right. \\
    \cdot \left. \left. \left( 1 + \frac{X ^{m ^3}}{1 - X ^{m ^3}} \left( 1 + \ldots
    \right) \right) \right) \right) = 1 + \sum \limits_{k \ge 0} \prod
    \limits_{i = 0}^{k} \frac{X ^{m ^i}}{1 - X ^{m ^i}}
\end{multline}
Let us denote
\[
    \prod \limits_{i = 0}^{k} \frac{X ^{m ^i}}{1 - X ^{m ^i}} = \sum \limits
    _{n \ge 0} P_k \n X^n,
\]
where $k \ge 0$. Clearly for any $k \ge 0$ the function $P_k : \Natural \to
\Natural$ is the polynomial of $\left( k + 1 \right)$ degree, $P_k \n[0] = 0$,
and the value $P_k \n$, $n \ge 0$, is equal to the number of partitions of $n$
into $\left( k + 1 \right)$ first powers of $m$, i.e.
\begin{equation} \label{eq:second_difference_P_k}
    P_k \n = \left| \left\{
        \begin{array}{*{20}c}
           {p_0, p_1, \ldots, p_k} & \vline & {\sum \limits_{i = 0}
           ^{k} p_i m ^i = n, p_i \ge 1, i = 0, 1, \ldots, k} \\
        \end{array}
    \right\} \right|
\end{equation}
It follows from~\eqref{eq:second_difference_partitions} and the definition of
$P_k$ that the following equalities hold
\begin{multline*}
    \Gamma ^{\n[2]} \n[X] = 1 + \sum \limits_{k \ge 0} \prod \limits_{i =
    0}^{k} \frac{X ^{m ^i}}{1 - X ^{m ^i}} = 1 + \sum \limits_{k \ge 0} \left(
    \sum \limits_{n \ge 0} P_k \n X ^n \right)\\
    = 1 + \sum \limits_{n \ge 1} \left( \sum \limits_{k \ge 0} P_k \n \right)
    X ^n = 1 + \sum \limits_{n \ge 1} \Growth{\Intermediate} ^{\n[2]} \n X ^n,
\end{multline*}
whence for all $n \ge 1$ we have
\[
    \Growth{\Intermediate} ^{\n[2]} \n = \sum \limits_{k \ge 0} P_k \n.
\]
It follows from~\eqref{eq:second_difference_P_k}, that the value
$\Growth{\Intermediate} ^{\n[2]} \n$ is equal to the number of partitions of
$n$ into ``sequential'' powers of $m$, that was required to be proved.

\par Corollary~\ref{cor:numerical_properties} is completely proved.
\end{proof}

\subsection{Asymptotics} \label{subsect:growth_asymptotics}

\par We quote the following result by Mahler \cite{Mahler1940}:
\begin{theorem} \label{th:special_equation}
Let $f \n[z]$ be a real function of the real variable $z \ge 0$ which in every
finite interval is bounded, but not necessarily continuous, and which satisfies
the equation
\[
    \frac{f \n[z + \omega] - f \n[z]}{\omega} = f \n[qz].
\]
If, as $z \to \infty$, $n$ is the integer for which
\[
    q ^{- \Argument{n - 1}} n \le z < q ^{- n} \Argument{n + 1},
\]
then
\[
    f \n[z] = \mathcal{O} \left( \frac{q^{\frac{1}{2} n \Argument{n - 1}}
    z^n}{n!} \right).
\]
This inequality can be improved to
\[
    f \n[z] = \frac{q^{\frac{1}{2} n \Argument{n - 1}}
    z^n}{n!} e ^{o \n[1] }.
\]
if $f \n[z]$ is greater than a positive constant $C$ for all sufficiently large
$z \ge 0$.
\end{theorem}

\begin{proof}[Proof of Theorem~\ref{th:estimates}]
Let us consider the function $f : \mathbb{R} ^{+} \to \mathbb{R}$, defined in
the following way
\begin{equation*}
    f \n[z] = \WordGrowthSemigroup{\Intermediate} \n[ {\left[ z \right]} ].
\end{equation*}
It follows from Item~\ref{cor_item:functional_equation} of
Corollary~\ref{cor:numerical_properties} that $f$ satisfies the conditions of
Theorem~\ref{th:special_equation} for $q = \frac{1}{m}$, $\omega = 1$. It
implies, that for all sufficiently large $l$ the equality holds
\[
    \WordGrowthSemigroup{\Intermediate} \n[l] = \frac{m ^{-\frac{1}{2} n
    \Argument{n - 1}} l^n}{n!} e ^{o \n[1] },
\]
where $n$ is defined by the inequalities $m ^{\Argument{n - 1}} n \le l < m
^{n} \Argument{n + 1}$.

\par It follows from~\eqref{eq:first_difference_aymptotics} (see
\cite{Mahler1940}, and also \cite{Pennington1953}) that logarithm of the word
growth function admits the following asymptotics
\begin{equation*}
    \log \WordGrowthSemigroup{\Intermediate} \n \sim \frac{\Argument{\log n
    }^2}{2 \log m},
\end{equation*}
whence
\[
    \WordGrowthSemigroup{\Intermediate} \n \sim \exp \left(
    \frac{\Argument{\log n} ^2}{2 \log m} \right) \sim n ^{\frac{\log n}{2
    \log m}}.
\]
It is proved in Corollary~\ref{cor:numerical_properties} that the equality
$\Growth{\Intermediate} \n = \left( 1/m \right)
\WordGrowthSemigroup{\Intermediate} \n[m \Argument{n + 1}] - 1/m$ holds for all
$n \ge 0$. Thus we have the sharp estimate
\[
    \Growth{\Intermediate} \n = \frac{1}{m} \WordGrowthSemigroup{\Intermediate}
    \n[m \Argument{n + 1}] - \frac{1}{m} \sim \frac{1}{m} \left( m \Argument{n
    + 1} \right) ^{\frac{\log \left( m \Argument{n + 1} \right) }{2 \log m}},
\]
with the ratios of left- to right-hand side tending to $1$ as $n \to \infty$.
As the functions $\Growth{\Intermediate}$ and $\GrowthSemigroup{\Intermediate}$
coincide then Theorem~\ref{th:estimates} is completely proved.
\end{proof}

\begin{proof}[Proof of Corollary~\ref{cor:growth_orders}]
It follows from Theorem~\ref{th:estimates} that the equality holds
\[
    \Growth{\Intermediate} \n = \GrowthSemigroup{\Intermediate} \n \sim
    \frac{1}{m} \left( m \Argument{n + 1} \right) ^{\frac{\log \left( m
    \Argument{n + 1} \right) }{2 \log m}},
\]
whence
\[
    \GrowthOrder{\Growth{\Intermediate}} = \GrowthOrder{
    \GrowthSemigroup{\Intermediate}} = \GrowthOrder{\frac{1}{m} \left( m
    \Argument{n + 1} \right) ^{\frac{\log \left( m \Argument{n + 1} \right) }{2
    \log m}}}.
\]
Two functions of a natural argument
\begin{align*}
    \Growth{1} \n & = n ^{\frac{\log n}{2 \log m}}, & \Growth{2} \n & =
    \frac{1}{m} \left( m \Argument{n + 1} \right) ^{\frac{\log \left( m
    \Argument{n + 1} \right) }{2 \log m}}.
\end{align*}
have the same growth orders, because they fulfilled the requirements of
Proposition~\ref{prop:the_same_growth_orders} for $h = \frac{1}{m}$, $a = m$,
$b = m$, $c = 0$. Therefore the equalities hold
\[
    \GrowthOrder{\Growth{\Intermediate}} =
    \GrowthOrder{\GrowthSemigroup{\Intermediate}} = \GrowthOrder{\frac{1}{m}
    \left( m \Argument{n + 1} \right) ^{\frac{\log \left( m \Argument{n + 1}
    \right) }{2 \log m}}} = \GrowthOrder{ n ^{\frac{\log n}{2 \log m}}},
\]
and the statement of the corollary is true.
\end{proof}

\section{The properties of $\SequenceOParam{\Intermediate}{2}$} \label{sect:sequence}

\begin{proof}[Proof of Theorem~\ref{th:sequence_growth}]
\ref{th_item:sequence_growth_orders}. It follows from
Corollary~\ref{cor:growth_orders} that for all $m \ge 2$ the equality holds
\[
    \GrowthOrder{\Growth{\Intermediate}} = \GrowthOrder{ n
    ^{\frac{\log n}{2 \log m}}},
\]
and therefore it's enough to prove that the following inequality holds
\begin{equation} \label{eq:growth_orders_inequality}
    \GrowthOrder{ n
    ^{\frac{\log n}{2 \log m}}} > \GrowthOrder{ n
    ^{\frac{\log n}{2 \log \left( m + 1 \right)}}}, \quad m \ge 2.
\end{equation}

\par Let us assume by contrary, that there exist positive numbers $C_1, C_2,
N_0 \in \Natural$ such that
\begin{equation} \label{eq:growth_orders_assumption}
    n ^{\frac{\log n}{2 \log m}} \le C_1 \Argument{C_2 n} ^{\frac{\log
    \Argument{C_2 n}}{2 \log \left( m + 1 \right)}}
\end{equation}
for any $n \ge N_0$. The functions at the left- and right-hand side are
positively defined non-decreasing functions, and the
assumption~\eqref{eq:growth_orders_assumption} is true if and only if the
inequality
\begin{equation*}
    \log \left( n ^{\frac{\log n}{2 \log m}} \right) \le \log \left( C_1
    \Argument{C_2 n} ^{\frac{\log \Argument{C_2 n}}{2 \log \left( m + 1
    \right)}} \right)
\end{equation*}
holds. The left-hand side can be transformed in the following way:
\begin{multline*}
    \log \left( C_1 \Argument{C_2 n} ^{\frac{\log \Argument{C_2 n}}{2 \log
    \left( m + 1 \right)}} \right) = \log C_1 + \frac{\log ^2 \Argument{C_2
    n}}{2 \log \left( m + 1 \right)} = \\
    = \log C_1 + \frac{1}{2 \log \left( m + 1 \right)} \left( \log ^2 n + 2
    \log C_2 \log n + \log ^2 C_2 \right).
\end{multline*}
Hence the assumption~\eqref{eq:growth_orders_assumption} is true if and only if
the following inequality holds
\begin{multline}  \label{eq:growth_orders_assumption2}
    \log ^2 n \left( \frac{1}{2 \log m} - \frac{1}{2 \log \left( m + 1 \right)}
    \right) - \log n \left( \frac{\log C_2}{\log \left( m + 1 \right)} \right) \\
    - \left( \log C_1 + \log ^2 C_2 \right) \le 0
\end{multline}
for all $n \ge N_0$. As $m \ge 2$ then the coefficient at $\log ^2 n$ satisfies
the inequality
\[
    \frac{1}{2 \log m} - \frac{1}{2 \log \left( m + 1 \right)} \ge 0,
\]
and therefore there exists $N_1 \in \Natural$, $N_1 \ge N_0$, such that the
inequality~\eqref{eq:growth_orders_assumption2} is false for all $n \ge N_1$.
Thus, we obtain the contradiction with the
assumption~\eqref{eq:growth_orders_assumption}, whence the
inequality~\eqref{eq:growth_orders_inequality} is true.
Item~\ref{th_item:sequence_growth_orders} of Theorem~\ref{th:sequence_growth}
is proved.

\par \ref{th_item:sequence_growth_functions}.
Furthermore, we separate the defining relations of different semigroups
$\Semigroup{\Intermediate}$ by the upper index $\n[m]$. Let us consider the set
of relations
\begin{equation} \label{eq:pointwise_defining_relations}
    \left\{ f_0 f_1 ^p f_0 = f_0, \, p \ge 0 \right\}.
\end{equation}
For fixed $p \ge 0$ the relation $f_0 f_1 ^p f_0 = f_0$ is the relation $R_A
\Pair{0}{p + 1}$, and it holds in each semigroup
$\Semigroup{\Intermediate[m]}$, where $m \ge {p + 2}$. On the other hand, the
defining relations of $\Semigroup{\Intermediate[m]}$ that don't belong to the
set~\eqref{eq:pointwise_defining_relations} can be applied to semigroup words
of length greater than $\Argument{m + 2}$. Therefore the
set~\eqref{eq:pointwise_defining_relations} can be considered as ``a pointwise
limit'' of the sets of defining relations
\begin{equation*}
    \left\{ R_A ^{\n[m]} \Pair{k}{p}, \, R_B ^{\n[m]} \Argument{k} \right\},
\end{equation*}
where $m \ge 2$, $k \ge 0$, $p = 1, 2, \ldots, {m - 1}$, as $m$ tends to
$+\infty$.

\par Let us consider the infinitely presented monoid
\[
    \Semigroup{} = \ATMonoid{ f_0 f_1 ^p f_0 = f_0, \, p \ge 0 },
\]
and we calculate its growth function $\Growth{\Semigroup{}}$. It follows from
the speculations above that $\Growth{\Semigroup{}}$ is the pointwise limit for
the growth function sequence
$\Sequence{\Growth{\Semigroup{\Intermediate}}}{2}$. It is easy to check that an
arbitrary element $\s \in \Semigroup{}$ can be unambiguously reduced to one of
the following forms
\begin{align*}
    & f_1 ^{p _0}, & p_0 & \ge 0,\\
    & f_1 ^{p _0} f_0 f_1 ^{p _1}, & p_0, p_1 & \ge 0.
\end{align*}
The word growth function $\WordGrowthSemigroup{}$ is defined by the following
equality
\[
    \WordGrowthSemigroup{} \n = \underbrace{1} _{f_1 ^n} + \underbrace{n}
    _{f_1 ^{p _0} f_0 f_1 ^{n - 1 - p _0}} = n + 1,
\]
whence
\[
    \Growth{\Semigroup{}} \n = \sum \limits _{i = 0} ^{n} \left( n + 1 \right)
    = \frac{\left( n + 1 \right) \left( n + 2 \right)}{2}.
\]
As $\Semigroup{\Intermediate}$ is a monoid for all $m \ge 2$ then the sequence
$\Sequence{\Growth{\Intermediate}}{2}$ tends pointwisely to the growth function
$\GrowthSemigroup{}$ as $m \to + \infty$, that is equal to $\frac{\left( n + 1
\right) \left( n + 2 \right)}{2}$.

\item \ref{th_item:sequence_automata}. Let $\xi$ be a cyclic permutation of
$\Alphabet{m}$ and $\theta$ be an identical permutation. Applying these
permutation to $\Intermediate$, we obtain the similar automaton
$\Intermediate'$ such that its automaton transformations have the following
decompositions
\begin{align*}
    f_0 & = ( f_0, e, e, \ldots, e ) \alpha_1, & f_1 & = ( f_1, e, e, \ldots, e
    ) \sigma.
\end{align*}
A pointwise limit of the automaton sequence $\Sequence{\Intermediate'}{2}$ is
the automaton $\Intermediate[\infty]'$ with the following automaton
transformations
\begin{align*}
    f_0 & = ( f_0, e, e, \ldots ) \left( {\begin{array}{*{20}c}
        {x_0} & {x_1} &  {x_2} & \ldots \\
        {x_1} & {x_1} &  {x_1} & \ldots \\
    \end{array} } \right), \\
    f_1 & = ( f_1, e, e, \ldots ) \left( {\begin{array}{*{20}c}
        {x_0} & {x_1} & {x_2} & \ldots \\
        {x_1} & {x_2} & {x_3} & \ldots \\
    \end{array} } \right).
\end{align*}
Moreover, it's convenient to consider the infinite alphabet $\Alphabet{}' = \{
x_{-1}, x_0, x_1, \ldots \}$, and we set up a bijection between $\Alphabet{}'$
and $\Alphabet{\infty} = \{x_0, x_1, x_2 \ldots \}$ in a natural way. Let
$\Intermediate[]'$ be an automaton shown on
Fig.\ref{fig:automaton_intermediate_infty}. It acts over the alphabet
$\Alphabet{}'$, and $\Intermediate[]'$ is a similar automaton to a pointwise
limit of $\Sequence{\Intermediate}{2}$.

Let $\Semigroup{\Intermediate[]'}$ be the automaton transformation monoid
defined by $\Intermediate[]'$. It is easily to check that the following
relations hold in $\Semigroup{\Intermediate[]'}$:
\begin{align*}
    f_0 f_1 ^{p} f_0 & = f_0, & p& \ge 0, && \text{and} & f_0 f_1 ^{p} & = f_0
    f_1, & p& \ge 1.
\end{align*}
Elements $f_1 ^p$, $f_1 ^p f_0$ and $f_1 ^p f_0 f_1$, $p \ge 0$ define pairwise
different automaton transformations over the set of infinite words over the
alphabet $\Alphabet{}'$. Thus $\Semigroup{\Intermediate[]'}$ has the following
presentation:
\[
    \Semigroup{\Intermediate[]'} = \ATMonoid{ f_0 f_1 ^p f_0 =
    f_0, p \ge 0, \, f_0 f_1 ^p = f_0 f_1, p \ge 1}.
\]
It follows from Item~\ref{th_item:sequence_growth_functions} that the monoid
$\Semigroup{}'$ is a factor-semigroup of the monoid $\Semigroup{}$ that can be
considered as a pointwise limit of the semigroup sequence
$\Sequence{\Semigroup{\Intermediate}}{2}$, but they are supposed to be
isomorphic. Moreover, for $n \ge 1$ there are $3n$ semigroup elements of length
$n$:
\begin{align*}
    & f_1 ^p, & p & = 0, 1, \ldots, n,\\
    & f_1 ^p f_0, & p & = 0, 1, \ldots, n - 1,\\
    & f_1 ^p f_0 f_1, & p & = 0, 1, \ldots, n - 2;
\end{align*}
whence the equality $\Growth{\Semigroup{}'} \n = 3n$ holds for all $n \ge 1$.
Obviously the growth functions $\Growth{\Semigroup{\Intermediate[]'}}$ and
$\Growth{\Semigroup{}}$ have different polynomial growth orders.

The theorem is completely proved.
\end{proof}

\section{Final remarks} \label{sect:final_remarks}

In the paper the sequence of the Mealy automata $\Intermediate$ is described.
From our point of view one of the most interesting properties is the property
of the growth function $\Growth{\Intermediate}$, $m \ge 2$, that is described
in Corollary~\ref{cor:numerical_properties},
Item~\ref{cor_item:functional_equation}:
\begin{equation*}
    \WordGrowthSemigroup{\Intermediate} \n[n + 1] =
    \WordGrowthSemigroup{\Intermediate} \n +
    \WordGrowthSemigroup{\Intermediate} \n[\Dividerm{n}],
\end{equation*}
where $n \ge 0$. Due to this fact, the second difference of the function
$\Growth{\Intermediate}$ is defined by the following equality
\[
    \Growth{\Intermediate} ^{\n[2]} \n = \WordGrowthSemigroup{\Intermediate} \n
    - \WordGrowthSemigroup{\Intermediate} \n[n - 1] =
    \WordGrowthSemigroup{\Intermediate} \n[\Dividerm{n - 1}]
\]
where $n \ge 1$, and, hence, we have
\[
    \Growth{\Intermediate} ^{\n[2]} \n[mn + 1] = \Growth{\Intermediate}
    ^{\n[2]} \n[mn + 2] = \ldots = \Growth{\Intermediate} ^{\n[2]} \n[mn + m]
\]
for all $n \ge 0$. Hence, the function $\Growth{\Intermediate} ^{\n[2]}$
consists of $m$ times repeated values of $\WordGrowthSemigroup{\Intermediate}$.

\begin{figure}[t]
  \centering
  \includegraphics*{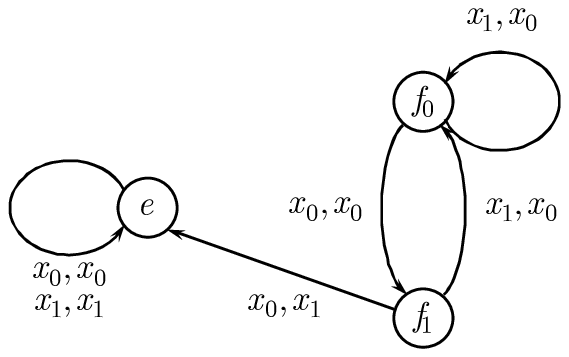}
  \caption{The automaton $\Automaton$}
  \label{fig:automaton_square}
\end{figure}

\par Let $\Automaton$ be the $3$-state Mealy automaton over the $2$-symbol
alphabet such that its Moore diagram is shown on
Figure~\ref{fig:automaton_square}. Let us denote its growth function by the
symbol $\Growth{\Automaton}$. The proposition holds

\begin{proposition}
The second difference $\Growth{\Automaton} ^{\n[2]}$ satisfies the following
equality
\[
    \Growth{\Automaton} ^{\n[2]} \n = 2 \Growth{\Automaton} ^{\n[2]} \n[n - 1]
    - \Growth{\Automaton} ^{\n[2]} \n[n - 2] + \Growth{\Automaton} ^{\n[2]}
    \n[ {\left[ \frac{n - 3}{2} \right]} ],
\]
where $n \ge 5$, and $\Growth{\Automaton} ^{\n[2]} \n[1] = 1$,
$\Growth{\Automaton} ^{\n[2]} \n[2] = 2$, $\Growth{\Automaton} ^{\n[2]} \n[3] =
3$, $\Growth{\Automaton} ^{\n[2]} \n[4] = 5$.
\end{proposition}

\par It follows from this proposition that the second difference of the
function $\Growth{\Automaton} ^{\n[2]}$, i.e. the fourth difference
$\Growth{\Automaton} ^{\n[4]}$, for $n \ge 5$ consists of doubled values of
$\Growth{\Automaton} ^{\n[2]}$. It is possible to assume that there exist
$3$-state Mealy automata such that the fourth difference of the growth function
consists of the second difference values that repeats $m$ times, where $m = 3,
4, \ldots$. Moreover, we put up the following problem:

\begin{problem}
Let $m \ge 2, k \ge 1$ be arbitrary positive integers. Do there exist Mealy
automata such that the $2k$-th finite difference of the growth function
consists of $m$ times repeated values of the $k$-th finite difference of this
growth function?
\end{problem}

\par It requires additional researches, but we think that the studying of Mealy
automata through the arithmetic properties of their growth functions and their
finite differences can produce many interesting examples.

\end{document}